\documentclass[12pt]{amsart}
\usepackage[all]{xy}
\usepackage{graphicx}
\usepackage{tikz}  

\usepackage{amsmath}
\usepackage{amssymb}
\usepackage{amsfonts}
\usepackage{mathrsfs}
\usepackage{mathtools}

\newcommand{\Ndb}{\ensuremath{\mathbb{N}}}

\newcommand{\norm}[1]{\Vert#1\Vert}

\newtheorem{theorem}{Theorem}[section]
\newtheorem{lemma}[theorem]{Lemma}
\newtheorem{corollary}[theorem]{Corollary}
\newtheorem{proposition}[theorem]{Proposition}
\newtheorem{definition}[theorem]{Definition}
\theoremstyle{remark}
\newtheorem{remark}[theorem]{\bf Remark}
\theoremstyle{definition}

\numberwithin{equation}{section}

\begin{document}

\title[Square functions for commuting families of Ritt operators]
{Square functions for commuting families of Ritt operators}

\author[O. Arrigoni]{Olivier Arrigoni}
\email{olivier.arrigoni@univ-fcomte.fr}
\address{Laboratoire de Math\'ematiques de Besan\c con, UMR 6623, 
CNRS, Universit\'e Bourgogne Franche-Comt\'e,
25030 Besan\c{c}on Cedex, FRANCE}

\date{\today}

\maketitle

\begin{abstract}

In this paper, we investigate the role of square functions defined for a $d$-tuple of commuting Ritt operators $(T_1,...,T_d)$ acting on a general Banach space $X$. Firstly, we prove that if the $d$-tuple admits a $H^\infty$ joint functional calculus, then it verifies various square function estimates. Then we study the converse when every $T_k$ is a $R$-Ritt operator. Under this last hypothesis, and when $X$ is a $K$-convex space, we show that square function estimates yield dilation of $(T_1,...,T_d)$ on some Bochner space $L_p(\Omega;X)$ into a $d$-tuple of isomorphisms with a $C(\mathbb{T}^d)$ bounded calculus. Finally, we compare for a $d$-tuple of Ritt operators its $H^\infty$ joint functional calculus with its dilation into a $d$-tuple of polynomially bounded isomorphisms.   
\end{abstract}


\vskip 1cm

\section{Introduction}

Many results on functional calculus are related to dilation results. Famous examples are given by Von Neumann's and Ando's inequalities. The latter uses in particular a joint dilation of two commuting contractions on a Hilbert space $H$ into a pair of commuting unitary operators on another Hilbert space $K$ which contains $H$ (see \cite{An}, \cite{P}). In \cite{ArLM}, the authors generalised this result provided that one considers a $d$-tuple of commuting Ritt contractions acting on a Hilbert space, $d \geq 2$. 

A fundamental paper of Fr\"ohlich and Weis (\cite{Fro}) shows that $H^\infty$ functional calculus implies dilation results for a sectorial operator $A$ acting on a Banach space $X$ having certain geometric properties. Under these conditions, they showed that an analytic semigroup whose negative generator $A$ has an $H^\infty$ functional calculus admits a dilation into a semigroup of isometries. 

Analogous results for a Ritt operator are proved in \cite{AFLM} and \cite{ArLM}. On one hand, \cite{AFLM} generalises the results of \cite{Fro}. On the other hand, \cite{ArLM} uses square function estimates related to each of the $T_k$'s and combines it to obtain joint dilation on Bochner spaces. 

In this paper, we pursue this study with square functions related to dilations for $d$-tuple of commuting Ritt operators, whereas the preceding results we cited use square functions for each of the commuting Ritt operators. As a continuation of the work of \cite{ArLM}, the main purpose of this paper is to study the relationship between $H^\infty$ functional calculus and dilations. We will obtain the following result.

\begin{theorem} \label{Th1}
Let $X$ be a reflexive $K$-convex Banach space and $p$ in $(1,\infty)$. Let $T=(T_1,...,T_d)$ be a $d$-tuple of commuting Ritt operators on $X$. Suppose that $T$ admits an $H^\infty(B_{\gamma_1} \times \cdots \times B_{\gamma_d})$ joint functional calculus for some $\gamma_1,...,\gamma_d$ in $(0,\frac{\pi}{2})$.

Then there exists a measure space $\Sigma$, a $d$-tuple of commuting isomorphisms $(U_1,...,U_d)$ on $L_p(\Sigma;X)$ admitting a $C(\mathbb{T}^d)$ bounded calculus and two bounded operators $J : X \to L_p(\Sigma;X)$ and $Q : L_p(\Sigma;X) \to X$ such that
\begin{equation} \label{C1}
T_1^{n_1} \cdots T_d^{n_d} = Q U_1^{n_1} \cdots U_d^{n_d} J, \qquad (n_1,...,n_d) \in \mathbb{N}^d.
\end{equation}
\end{theorem}

Next we study a form of converse, that is a dilation property implies $H^\infty$ functional calculus provided we consider $R$-Ritt operators.

\begin{theorem} \label{Th2}

Let $X$ be a Banach space and $p \in (1,\infty)$. Let $(T_1,...,T_d)$ be a $d$-tuple of commuting operators acting on $X$ such that every $T_k$ is an $R$-Ritt operator, $k=1,...,d$. Suppose that there exist a measure space $\Sigma$, a $d$-tuple of commuting isomorphisms $(U_1,...,U_d)$ acting on $L_p(\Sigma;X)$ having a $C(\mathbb{T}^d)$ bounded calculus and two bounded operator $J : X \to L_p(\Sigma;X)$, $Q : L_p(\Sigma;X) \to X$ such that (\ref{C1}) is verified. 

Then  there exist $\gamma_1,...,\gamma_d$ in $(0,\frac{\pi}{2})$ such that $(T_1,...,T_d)$ admits an $H^\infty(B_{\gamma_1} \times \cdots \times B_{\gamma_d})$ joint functional calculus.

\end{theorem}

Recall now that according to \cite[Proposition 7.7]{LM}, any polynomially bounded $R$-Ritt operator admits an $H^\infty$ functional calculus . As a key step of the proof of Theorem \ref{Th1}, we establish a multivariate version for $d$-tuples as follows.

\begin{theorem} \label{Th3}

Let $X$ be a Banach space and $(T_1,...,T_d)$ a $d$-tuple of commuting operators acting on $X$ such that every $T_k$ is an $R$-Ritt operator, $k=1,...,d$. Suppose that $(T_1,...,T_d)$ is polynomially bounded, that is there exists a constant $C \geq 1$ such that for any polynomial function $h$ of $d$ variables we have
\begin{equation*} 
\left\| h(T_1,...,T_d) \right\| \leq C \text{sup} \left\lbrace \left| h(z_1,...,z_d) \right| : (z_1,...,z_d) \in \mathbb{T}^d \right\rbrace.
\end{equation*}

Then  there exist $\gamma_1,...,\gamma_d$ in $(0,\frac{\pi}{2})$ such that $(T_1,...,T_d)$ admits a $H^\infty(B_{\gamma_1} \times \cdots \times B_{\gamma_d})$ joint functional calculus.
\end{theorem}

We mention that in the recent paper  \cite{ArLM2}, tools and results similar to the ones presented here have been introduced and developed for $d$-tuples of commuting semigroups and sectorial operators.


We now give a brief description of this paper. In Section \ref{SF}, we give background on $H^\infty$ joint functional calculus introduced in \cite{ArLM} and we define square functions on a general Banach space related to a tuple of commuting Ritt operators. We make them explicit for Hilbert spaces and Banach lattices (and especially for the $L_p$-spaces). Section \ref{JointSquare} shows that any tuple of commuting Ritt operators having a $H^\infty$ joint functional calculus verifies a square function estimate using the quadratic functional calculus. In Section \ref{SquareJoint}, we study the converse result, that is for which case square function estimates imply the $H^\infty$ joint functional calculus property. In Section \ref{Dilations}, we apply square function estimates and prove Theorem \ref{Th1} above.  Section \ref{HRRitt} is devoted to Theorems \ref{Th2} and \ref{Th3}. It uses an extension of a result of Franks and McIntosh. This latter provides a way to reduce the domain of $H^\infty$ functional calculus for operators acting on an Hilbert space.


We end this section by fixing some notations. Throughout 
we let $B(X)$ denote the Banach algebra of 
all bounded operators on some Banach space $X$. We let $I_X$ 
denote the identity operator on $X$. 
For any (possibly unbounded) operator $A$ on $X$, we let 
$\sigma(A)$ denote the spectrum of $A$ and for every $\lambda$ in $\mathbb{C} 
\setminus \sigma(A)$, we let  $R(\lambda,A)=(\lambda I_X - A)^{-1}$ denote the
resolvent operator. Next, we let $\text{Ker}(A)$ and 
$\text{Ran}(A)$ denote the kernel and the range of $A$, respectively.

For any $a \in \mathbb{C}$ and $r>0$, $D(a,r)$ will denote the open disc 
centered at $a$ with radius $r$. Then we let $\mathbb{D} = D(0,1)$ denote
the unit disc of $\mathbb{C}$ and we set $\mathbb{T} = \overline{\mathbb{D}} 
\setminus \mathbb{D}$. 

If $X$ is a Banach space and $\mathcal{O}$ is an open non empty subset of $\mathbb{C}^d$, for some integer 
$d \geq 1$, we will denote by $H^\infty(\mathcal{O};X)$ the Banach space of all bounded 
holomorphic functions $f \colon \mathcal{O} \to X$, equipped with the norm 
\begin{equation*}
\left\| f \right\|_{\infty,\mathcal{O}} = \text{sup} \left\lbrace \left\| 
f(z_1,\ldots,z_d) \right\| : (z_1,\ldots,z_d) \in \mathcal{O} \right\rbrace.
\end{equation*}  
When $X = \mathbb{C}$, we will simply write $H^\infty(\mathcal{O};\mathbb{C}) = H^\infty(\mathcal{O})$ which is a Banach algebra for the pointwise multiplication and the preceding norm.

If $X$ is a Banach space, $(\Omega,\mu)$ is a measure space and $p \in(1,\infty)$, 
we denote by $L_p(\Omega;X)$ the Bochner space of all classes of measurable functions 
$f \colon \Omega \to X$ such that $\int_\Omega \left\| 
f(\omega) \right\|^p d\mu(\omega) < \infty$, and we let 
$L_p(\Omega) = L_p(\Omega;\mathbb{C}) $. We refer the reader e.g. to \cite{Hyt1} 
for more details.

For any finite set $\Lambda$, we will denote by $|\Lambda|$ the number of elements 
of $\Lambda$ and if $\Lambda'$ is a susbset of $\Lambda$, we usually denote by $(\Lambda')^c = \Lambda \setminus \Lambda'$ the complementary of $\Lambda'$ in $\Lambda$.

The set of nonnegative integers will be denoted by $\mathbb{N} = 
\left\lbrace 0,1,2,... \right\rbrace$. We set $\Ndb^*=\Ndb\setminus\{0\}$.

In certain proofs, we use the notation $\lesssim$ to indicate an 
inequality valid up to a constant which does not depend on the particular 
elements to which it applies. 
We use as well notation $A \simeq B$ to say that we have both $A \lesssim B$ and $B \lesssim A$.

\section{Square functions on general Banach spaces} \label{SF}

In this section, we introduce the square functions related to a commuting family of Ritt operators on general Banach spaces. We also recall basic definitions and properties of the $H^\infty$ joint functional calculus. For proofs and details on this subject, we refer the reader to \cite[Section 2]{ArLM}.

A bounded operator $T \colon X \to X$ is a Ritt operator provided there exists a constant $C>0$ such that
\begin{equation*}
\left\| T^n \right\| \leq C \qquad\hbox{and}\qquad 
\left\| n(T^n-T^{n-1}) \right\| \leq C, \qquad n \geq 1.
\end{equation*}

Ritt operators have a spectral characterisation. 
Namely $T$ is a Ritt operator if and only if $\sigma(T) \subset \overline{\mathbb{D}}$ 
and there exists a constant $K>0$ such that
\begin{equation*}
\left\| (\lambda-1)R(\lambda,T) \right\| \leq K, \qquad \lambda \in \mathbb{C},\ |\lambda| > 1.
\end{equation*}

For any  $a$ in $(0,\frac{\pi}{2})$, let $B_a$ denote the \textit{Stolz domain} of angle $a$, 
defined as the interior of the convex hull of $1$ and the disc $ D(0,\sin(a))$.

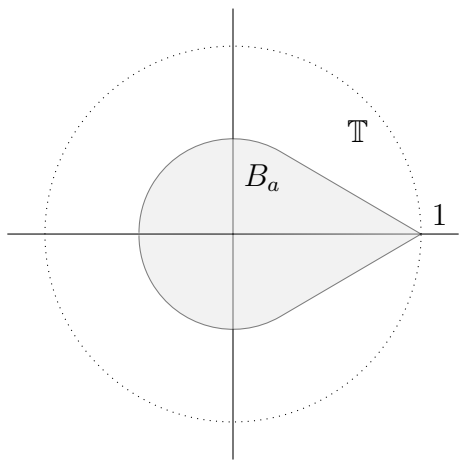
\begin{center}
\begin{tikzpicture}
[scale = 2.5];
\draw (-1.2,0)--(1.2,0);
\draw (0,-1.2)--(0,1.2);
\draw [dotted] circle(1);

\draw [fill=gray!20,opacity=0.5] (1,0)--(0.25,0.44)--(0.25,0.44) arc (60:180:0.5);
\draw [fill=gray!20,opacity=0.5] (1,0)--(0.25,-0.44)--(0.25,-0.44) arc (-60:-180:0.5);

\draw [right] (0,0.3) node{$B_a$};

\draw [right] (0.55,0.55) node{$\mathbb{T}$};

\draw [right] (1,0.1) node{$1$};

\end{tikzpicture}
\end{center}

Note that for any Stolz domain $B_a$, the set  $\left\lbrace \frac{\left|1-z\right|}{1-\left|z\right|} : z \in B_a \setminus \left\lbrace 1 \right\rbrace \right\rbrace $ is bounded (see e.g \cite{GT} for the result, as well as complements).

It turns out that if $T$ is a Ritt operator, then $\sigma(T) \subset \overline{B_a}$ for some
$a$ in $(0,\frac{\pi}{2})$.
More precisely (see \cite[Lemma 2.1]{LM}), one can find $a \in (0,\frac{\pi}{2})$ 
such that $\sigma(T) \subset \overline{B_a} $ and for any $b$ in $(a,\frac{\pi}{2})$, 
there exists a constant $K_b>0$ such that
\begin{equation} \label{majRittbeta}
\left\| (\lambda-1)R(\lambda,T) \right\| \leq K_{b}, 
\qquad \lambda \in \mathbb{C} \setminus \overline{B_b}.
\end{equation}
If this property holds, then we say that $T$ is a Ritt operator of type $a$. 
We refer to \cite{Ly,NZ,Nev} for the facts above and also to \cite{LM}
and the references therein for complements on the class of Ritt operators.

We now define the joint functional calculus of a family of commuting Ritt operators.

Let $d \geq 1$ be an integer and 
let $\gamma_1,\ldots,\gamma_d$ be elements of $(0,\frac{\pi}{2})$. For any 
subset $\Lambda$ of $\left\lbrace 1,\ldots,d \right\rbrace $, we denote by $H^\infty_0 
\left( \prod_{i \in \Lambda} B_{\gamma_i} \right) $ the subalgebra of $H^\infty 
\left( B_{\gamma_1} \times \cdots \times B_{\gamma_d} \right) $ of all holomorphic bounded 
functions $f$ depending only on variables 
$(\lambda_i)_{i\in \Lambda}$ and such that there exist positive constants $c$ and $(s_i)_{i \in \Lambda}$ verifying
\begin{equation} \label{Hinfini0Ritt}
\left| f(\lambda_1,\ldots,\lambda_d) \right| \leq c ~ \prod_{i \in \Lambda} 
|1-\lambda_i|^{s_i}, \qquad (\lambda_i)_{i\in \Lambda} \in \prod_{i \in \Lambda} B_{\gamma_i}.
\end{equation}
When $\Lambda=\emptyset$, $H^\infty_0 \left( \prod_{i \in \emptyset} B_{\gamma_i} \right) $ is the space of constant 
functions on $B_{\gamma_1} \times \cdots \times B_{\gamma_d}$.

Let $(T_1,\ldots,T_d)$ be a $d$-tuple of commuting Ritt operators. Assume that for any $k=1,\ldots,d$,
$T_k$ is of type $a_k \in (0,\gamma_k)$ and let $b_k \in (a_k,\gamma_k)$.

For any $f$ in $H_0^{\infty}(\prod_{i \in \Lambda} B_{\gamma_i})$ with 
$\Lambda \subset \left\lbrace 1,\ldots,d \right\rbrace$, $\Lambda \neq \emptyset$, we let
\begin{equation} \label{fTi}
f(T_1,\ldots,T_d) = \left( \dfrac{1}{2\pi i} \right)^{|\Lambda|} \int_{\prod_{i \in \Lambda} \partial B_{b_i}} f(\lambda_1,\ldots,\lambda_d) \prod_{i \in \Lambda} R(\lambda_i,T_i) \prod_{i \in \Lambda} d \lambda_i,
\end{equation}
where the $\partial B_{b_i}$ are oriented counterclockwise for all $i$ in $\Lambda$. 
This integral is absolutely convergent, hence defines an element of $B(X)$, its definition does not depend on the $b_i$ and 
the linear mapping $f \mapsto f(T_1,\ldots,T_d)$ is an algebra homomorphism 
from $H_0^{\infty}(\prod_{i \in \Lambda} B_{\gamma_i})$ into $B(X)$. If $f \equiv c$ is a constant function, 
then we let $f(T_1,\ldots,T_d) = cI_X$.

Next we define
$$
H_{0,1}^\infty(B_{\gamma_1} \times \cdots \times B_{\gamma_d}) = 
\bigoplus_{\Lambda \subset \left\lbrace 1,\ldots,d \right\rbrace} 
H^\infty_0 \left( \prod_{i \in \Lambda} B_{\gamma_i} \right)
$$
where the sum above is indeed a direct one. 

For any function $f = \sum_{\Lambda \subset \left\lbrace 1,\ldots,d \right\rbrace} f_\Lambda$
in $H_{0,1}^\infty(B_{\gamma_1} \times \cdots \times B_{\gamma_d})$,
with $f_\Lambda \in H^\infty_0 \left( \prod_{i \in \Lambda} B_{\gamma_i} \right)$, 
we let $f(T_1,\ldots,T_d) = \sum_{\Lambda \subset \left\lbrace 1,\ldots,d \right\rbrace} f_\Lambda(T_1,\ldots,T_d)$,
where every $f_\Lambda(T_1,\ldots,T_d)$ is defined by (\ref{fTi}). The mapping
$f\mapsto f(T_1,\ldots,T_d)$ is called the functional calculus mapping associated to
$(T_1,\ldots,T_d)$. This is an algebra homomorphism from $H_{0,1}^\infty(B_{\gamma_1} \times \cdots \times B_{\gamma_d})$ 
into $B(X)$. 

\begin{definition}
We say that $(T_1,\ldots,T_d)$ admits an $H^\infty(B_{\gamma_1} \times \cdots \times B_{\gamma_d})$ 
joint functional calculus if the above functional calculus mapping is bounded,
that is, there exists a constant $K >0$ such that for every $f$ 
in $H_{0,1}^\infty(B_{\gamma_1} \times \cdots \times B_{\gamma_d})$, we have
\begin{equation} \label{calcjointRitt}
\left\| f(T_1,\ldots,T_d) \right\| \leq K \left\| f \right\|_{\infty, B_{\gamma_1} 
\times \cdots \times B_{\gamma_d}}.
\end{equation} 
\end{definition}

We observe that $(T_1,\ldots,T_d)$
admits an $H^\infty(B_{\gamma_1} \times \cdots \times B_{\gamma_d})$ 
joint functional calculus if and only if 
$f\mapsto f(T_1,\ldots, T_d)$ is bounded on 
$H^\infty_0 \left( \prod_{i \in \Lambda} B_{\gamma_i} \right)$ for any $\Lambda \subset\{1,\ldots,d\}$.

Further if $(T_1,\ldots,T_d)$ admits an $H^\infty(B_{\gamma_1} \times \cdots \times 
B_{\gamma_d})$  joint functional calculus, then 
for every $k=1,\ldots,d$, $T_k$ admits an $H^\infty(B_{\gamma_k})$ functional calculus 
in the  sense of \cite[Definition 2.4]{LM}. More generally, any subfamily of $(T_1,...,T_d)$ admits a joint functional calculus if $(T_1,...,T_d)$ admits one.

By \cite[Proposition 2.5]{ArLM}, it suffices to obtain inequality (\ref{calcjointRitt}) only for polynomial functions $f$ of $d$ variables to prove that $(T_1,\ldots,T_d)$ admits an $H^\infty(B_{\gamma_1} \times \cdots \times B_{\gamma_d})$ 
joint functional calculus.

In order to define square functions for commuting families of Ritt operators, we recall some background on Rademacher averages on a general Banach space $X$.
Let $\mathcal{I}$ be a nonempty countable set. Let $(r_{\iota})_{\iota \in \mathcal{I}}$ be a family of independent Rademacher variables indexed by $\mathcal{I}$ on some probability space $\Omega_0$. If $1 < p < \infty$, we denote by $\text{Rad}_p(\mathcal{I};X)$ the closed subspace of $L_p(\Omega_0;X)$ which is the closure of the linear span of all the finite sums of type 
\begin{equation*}
\sum_{\iota \in \mathcal{I}} r_{\iota} \otimes x_{\iota},
\end{equation*}
where $(x_{\iota})$ is a finite family of $X$.

According to the Khintchine Kahane inequality, all the spaces $\text{Rad}_p(\mathcal{I};X)$ are isomorphic for $1 < p < \infty$. We denote by $\text{Rad}(\mathcal{I};X)$ the space $\text{Rad}_2(\mathcal{I};X) \subset L_2(\Omega_0;X)$. In the case where $\mathcal{I} = \mathbb{N}^*$, we set $\text{Rad}_p(X) = \text{Rad}_p(\mathbb{N}^*;X)$ and $\text{Rad}(X) = \text{Rad}_2(X)$.

In particular, when $\mathcal{I} = (\mathbb{N}^*)^d$, we let $(r_{k_1,...,k_d})$ denote an independent family of Rademacher variables indexed by $(\mathbb{N}^*)^d$.

Theorem \ref{Kwapien} below gives a criterion of summability, which is a generalisation of Kwapie{\'n}'s Theorem saying that if $X$ is a Banach space which does not contain $c_0$, a series $\sum r_n \otimes x_n $ converges in Rad($X$) if and only if its partial sums are uniformely bounded.

This generalisation uses the following lemma (see \cite[Prop. 6.1.5]{Hyt2}).

\begin{lemma} \label{indsymetric}
Let $X$ be a Banach space and $(\Omega,\mathbb{P})$ be a measure space. Let $\xi$ and $\eta$ be random variables from $\Omega$ to $X$. If $\eta$ is real-symmetric (i.e $\eta$ and $-\eta$ are identically distributed) and independant of $\xi$, then for all $1 \leq p \leq \infty$ we have
\begin{equation} \label{}
\left\| \xi \right\|_{L^p(\Omega;X)} \leq \left\| \xi + \eta \right\|_{L^p(\Omega;X)}.
\end{equation}

\end{lemma}

\begin{theorem} \label{Kwapien}

Let $X$ be a Banach space which does not contain $c_0$. Let $d \geq 2$ be an integer and $(x_{k_1,...,k_d})_{(k_1,...,k_d) \in (\mathbb{N}*)^d}$ be a family of $X$. Suppose that there exists a constant $K \geq 0$ such that 
\begin{equation} \label{sommepartielle}
\left\| \sum_{1 \leq k_1,...,k_d \leq N} r_{k_1,...,k_d} \otimes x_{k_1,...,k_d} \right\|_{\text{Rad}((\mathbb{N}^*)^d;X)} \leq K, \qquad N \in \mathbb{N}^*.
\end{equation}

Then the family $(r_{k_1,...,k_d} \otimes x_{k_1,...,k_d})_{(k_1,...,k_d) \in (\mathbb{N}*)^d}$ is summable in $\text{Rad}((\mathbb{N}^*)^d;X) $.

\end{theorem}

\paragraph*{} \textit{Proof of Theorem \ref{Kwapien} :} Let $\psi : \mathbb{N}^* \to (\mathbb{N}^*)^d$ be an arbitrary bijection. We prove that $\sum r_{\psi(n)} \otimes x_{\psi(n)} $ converges in Rad($X$). By Kwapie{\'n}'s Theorem, it suffices to prove that the partial sums of this series are uniformely bounded. Let $M \geq 1$ be an integer. There exists $N \geq 1$ such that
\begin{equation*} \label{bijection}
\left\lbrace \psi(1),...,\psi(M) \right\rbrace \subset \left\lbrace 1,...,N \right\rbrace^d.
\end{equation*}

Then by Lemma \ref{indsymetric} and assumption (\ref{sommepartielle}), we have

\begin{equation*}
\left\| \sum_{n=1}^{M} r_{\psi(n)} \otimes x_{\psi(n)} \right\|_{\text{Rad}(X)} \leq \left\|  \sum_{1 \leq k_1,...,k_d \leq N} r_{k_1,...,k_d} \otimes x_{k_1,...,k_d} \right\|_{\text{Rad}((\mathbb{N}^*)^d;X)} \leq K.
\end{equation*}

This yields the result. $\square$

For $\alpha>0$, define $p_\alpha : z \mapsto (1-z)^\alpha$ which is an element of $H_0^\infty(B_\gamma)$ for every $\gamma$ in $(0,\frac{\pi}{2})$. Then we let $(I_X-T)^\alpha = p_\alpha(T)$ in the sense of the functional calculus defined in (\ref{fTi}). 

We now define the square functions for commuting Ritt operators. Let $T=(T_1,...,T_d)$ be a $d$-tuple of commuting Ritt operators on $X$ and $\alpha=(\alpha_1,...,\alpha_d)$ in $(\mathbb{R}_+^*)^d$. Let $\Lambda$ be a subset of $\left\lbrace 1,...,d \right\rbrace$. Define 
\begin{equation*}
\alpha_\Lambda = (\alpha_i)_{i \in \Lambda} \in (\mathbb{R}_+^*)^\Lambda.
\end{equation*}
Then we let $(r_{(k_i)_{i \in \Lambda}})$ be a family of independent Rademacher variables indexed by $(\mathbb{N}^*)^\Lambda$.
Define now for any $(k_i)_{i \in \Lambda} \in (\mathbb{N}^*)^\Lambda$ and $x$ in $X$
\begin{equation*}
x_{(k_i)} = \left[ \prod_{i \in \Lambda} k_i^{\alpha_i-\frac{1}{2}} T_i^{k_i-1} (I_X-T_i)^{\alpha_i} \right] x.
\end{equation*}

We let for any $x$ of $X$
\begin{equation} \label{fonccarreRittLambda}
\norm{x}_{T,\alpha_\Lambda} = \left\| \sum_{(k_i)_{i \in \Lambda} \in (\mathbb{N}^*)^\Lambda} r_{(k_i)} \otimes x_{(k_i)} \right\|_{\text{Rad}((\mathbb{N}^*)^\Lambda;X)},
\end{equation}
if the family $\left( r_{(k_i)} \otimes x_{(k_i)} \right)_{(k_i) \in (\mathbb{N}^*)^\Lambda}$ is summable 
and we let $\norm{x}_{T,\alpha_\Lambda} = \infty $ otherwise. 

If $\Lambda = \left\lbrace 1,...,d \right\rbrace$, we will simply use the notation
\begin{equation} \label{fonccarreRitt}
\norm{x}_{T,\alpha} = \left\| \sum_{k_1,...,k_d \geq 1} \left( \prod_{i=1}^d k_i^{\alpha_i-\frac{1}{2}} \right) r_{k_1,...,k_d} \otimes \left( \prod_{i=1}^d T_i^{k_i-1} (I_X-T_i)^{\alpha_i} x  \right)\right\|_{\text{Rad}((\mathbb{N}^*)^d;X)}.
\end{equation} 

In the case where $X = E(S)$ is a Banach lattice of functions with finite cotype, the Khintchine-Maurey inequalities imply that the family $(r_{(k_i)} \otimes x_{(k_i)})_{(k_i) \in (\mathbb{N}^*)^\Lambda}$ is summable if and only if $\left\| \left( \sum_{(k_i) \in (\mathbb{N}^*)^\Lambda} \left| x_{(k_i)} \right|^2 \right)^\frac{1}{2} \right\|_{E(S)}$ is finite and in this case, we have
\begin{equation*}
\norm{x}_{T,\alpha_\Lambda} \simeq \left\| \left( \sum_{(k_i) \in (\mathbb{N}^*)^\Lambda} \left| x_{(k_i)} \right|^2 \right)^\frac{1}{2} \right\|_{E(S)};
\end{equation*}
(see \cite{Hyt2}).

We are interested in the $d$-tuples $(T_1,...,T_d)$ for which there exists a constant $K >0$ such that for any subset $\Lambda$ of $\left\lbrace 1,...,d \right\rbrace $ and for any $x$ in $X$, we have
\begin{equation} \label{squareRitt}
\norm{x}_{T,\alpha_\Lambda} \leq K \norm{x},
\end{equation}
for some $d$-tuple $(\alpha_1,...,\alpha_d)$. If such an inequality (\ref{squareRitt}) happens, we will say that the $d$-tuple $(T_1,...,T_d)$ admits a \textit{square function estimate}.

Note that this square function estimate depends a priori on the $d$-tuple $\alpha$. We will get back to this problem with Theorem \ref{th21} below.

\section{From $H^\infty$ joint functional calculus to square function estimates} \label{JointSquare}

In this section, we aim to show the following theorem.

\begin{theorem} \label{Thsquarefunc}
Let $X$ be a Banach space with a finite cotype. 
Suppose that $T=(T_1,...,T_d)$ is a $d$-tuple of commuting Ritt operators on $X$ which has an $H^\infty(B_{\gamma_1} \times \cdots \times B_{\gamma_d})$ joint functional calculus for some $\gamma_1,...,\gamma_d $ in $(0,\frac{\pi}{2})$. Let $\alpha = (\alpha_1,...,\alpha_d)$ in $(\mathbb{R}_+^*)^d $. Then $T$ verifies the following square function estimate : there exists a constant $K>0$ such that for any $x$ in $X$ we have
\begin{equation} \label{squareestRitt}
\left\| x \right\|_{T,\alpha} \leq K \left\| x \right\|.
\end{equation}

\end{theorem}

\begin{remark} \label{Remarque}
We recall that every subfamily of a $d$-tuple $(T_1,...,T_d)$ as in Theorem \ref{Thsquarefunc} has an $H^\infty$ joint functional calculus too. Thus, one can obtain the inequality (\ref{squareestRitt}) replacing $\alpha$ by $\alpha_\Lambda$, $\Lambda \subset \left\lbrace 1,...,d \right\rbrace$, where we refer to (\ref{fonccarreRittLambda}) for the definition of square function.
\end{remark}

This result appeals to the notion of quadratic functional calculus which is defined as follows.

\begin{definition}

Let $X$ be a Banach space. Let $T=(T_1,...,T_d)$ be a $d$-tuple of commuting Ritt operators on $X$ such that $T_k$ is of type $a_k$ for $k=1,...,d$.
Let $\gamma_k \in (a_k,\frac{\pi}{2})$ for $k=1,...,d$. We say that $T$ admits a quadratic $H^\infty (B_{\gamma_1} \times \cdots \times B_{\gamma_d})$ functional calculus if there exists a constant $C>0$ such that for any finite family $(\varphi_{i})_{i \in \mathcal{I}}$ in $H^\infty_{0,1}(B_{\gamma_1} \times \cdots \times B_{\gamma_d})$ and $x$ in $X$,
\begin{equation} \label{estquadraRitt}
\left\| \sum_{i \in \mathcal{I}} r_{i} \otimes \varphi_{i}(T_1,...,T_d) (x) \right\|_{\textnormal{Rad}(\mathcal{I};X)} \leq C \left\| x \right\| \left\| \left( \sum_{i \in \mathcal{I}} \left| \varphi_{i} \right|^2 \right)^\frac{1}{2} \right\|_{\infty,B_{\gamma_1} \times \cdots \times B_{\gamma_d}},
\end{equation}
with $(r_i)_{i \in \mathcal{I}}$ a family of independent Rademacher variables.

\end{definition}

Theorem \ref{Thsquarefunc} is obtained by combining the following two propositions.

\begin{proposition} \label{Calcquadra}
Let $X$ be a Banach space with a finite cotype. Let $T=(T_1,...,T_d)$ be a $d$-tuple of commuting Ritt operators on $X$. Suppose that $T$ has an $H^\infty(B_{b_1} \times \cdots \times B_{b_d})$ joint functional calculus for some $b_1,...,b_d$ in $(0,\frac{\pi}{2})$. Then for any $\gamma_1,...,\gamma_d$ such that $\frac{\pi}{2} > \gamma_k > b_k$ for $k=1,...,d$, $T$ has a quadratic $H^\infty(B_{\gamma_1} \times \cdots \times B_{\gamma_d})$ functional calculus.
\end{proposition}

\begin{proposition} \label{Quadrasquare}
Let $X$ be a Banach which does not contain $c_0$. Let $T=(T_1,...,T_d)$ be a $d$-tuple of commuting Ritt operators on $X$. Suppose that $T$ has a quadratic $H^\infty(B_{\gamma_1} \times \cdots \times B_{\gamma_d})$ functional calculus for some $\gamma_1,...,\gamma_d$ in $(0,\frac{\pi}{2})$. Let $\alpha=(\alpha_1,...,\alpha_d)$ in $(\mathbb{R}_+^*)^d$. Then $T$ satisfies a square function estimate, that is there exists a constant $K>0$ such that
\begin{equation} \label{estcarre}
\left\| x \right\|_{T,\alpha} \leq K \left\| x \right\|.
\end{equation}

\end{proposition}

\paragraph*{} \textit{Proof of Proposition \ref{Calcquadra} : } The first step of the proof relies on a decomposition principle for holomorphic functions of several variables. The original idea of such a decomposition is due to Franks and McIntosh (see \cite{FMI}). We have a decomposition for Stolz domain (see \cite[Section 6]{ArLM}) which is useful for $d$-tuple of Ritt operators.

We let $0 < b_k < \gamma_k < \frac{\pi}{2}$ be angles for $k=1,...,d$.

The decomposition in \cite[Section 6]{ArLM} provides sequences of holomorphic functions $(\Psi_{k,i_k})_{i_k \in \mathbb{N}}$ and $(\tilde{\Psi}_{k,i_k})_{i_k \in \mathbb{N}}$ in $H_0^\infty(B_{b_k})$ for $k=1,...,d$ such that the two following properties hold.

\begin{itemize}
\item[(i)] For every $p>0$, there exists a constant $C_p>0$ such that 
\begin{equation} \label{majphiij1}
\text{sup} \left\lbrace \sum_{i_k = 1}^{\infty} \left| \Psi_{k,i_k}(\zeta_k)\right|^p : \zeta_k \in B_{b_k} \right\rbrace \leq C_p, \qquad k=1,...,d,
\end{equation}
\begin{equation} \label{majtildephiij}
\text{sup} \left\lbrace \sum_{i_k = 1}^{\infty} \left| \tilde{\Psi}_{k,i_k}(\zeta_k)\right|^p : \zeta_k \in B_{b_k} \right\rbrace \leq C_p, \qquad k=1,...,d;
\end{equation} 

\item[(ii)] For any Banach space $Z$ and a function $h$ in $H^\infty(B_{\gamma_1} \times \cdots \times B_{\gamma_d} ; Z) $, there exists a family $(a_{i_1,...,i_d})$ in $Z$ indexed by $\mathbb{N}^d$ such that for every $(\zeta_1,...,\zeta_d) $ in $B_{b_1} \times \cdots \times B_{b_d}$ we have
\begin{equation} \label{FMdvarialbes}
h(\zeta_1,\ldots,\zeta_d) = \sum_{i_1,\cdots,i_d} a_{i_1,...,i_d} \Psi_{1,i_1}(\zeta_1)\tilde{\Psi}_{1,i_1}(\zeta_1) \cdots \Psi_{d,i_d}(\zeta_d)\tilde{\Psi}_{d,i_d}(\zeta_d),
\end{equation}
and there exists a constant $C>0$ (independent of $h$) such that for every $(i_1,...,i_d) $ in $(\mathbb{N}^*)^d $ and $k=1,...,d$,
\begin{equation} \label{majalphaij1}
\left\| a_{i_1,...,i_d} \right\| \leq C \left\| h \right\|_{\infty, B_{b_1} \times \cdots \times B_{b_d}}.
\end{equation}
\end{itemize}

The original proof is done for $Z = \mathbb{C}$ in \cite[Section 6]{ArLM}. However, This proof works as well for any Banach space $Z$ and for any function $h$ in $H^\infty(B_{\gamma_1} \times \cdots \times B_{\gamma_d} ; Z) $.

Let now $(\eta_{i_1,...,i_d})$ be a finite family of complex numbers and let $m \geq 1$ be an integer. As $(T_1,...,T_d)$ has an $H^\infty(B_{b_1} \times \cdots \times B_{b_d})$ joint functional calculus, we have the estimate
\begin{align*}
\left\| \sum_{i_1,...,i_d=1}^m \eta_{i_1,...,i_d} \Psi_{1,i_1}(T_1)  \cdots \Psi_{d,i_d}(T_d) \right\| 
         & \lesssim \underset{(z_1,...,z_d) \in B_{b_1} \times \cdots \times B_{b_d} }{\text{sup}} \left| \sum_{i_1,...,i_d=1}^m \eta_{i_1,...,i_d} \Psi_{1,i_1}(z_1) \cdots \Psi_{d,i_d}(z_d) \right| \\ 
 		 &  \leq \text{sup} \left| \eta_{i_1,...,i_d} \right| \underset{B_{b_1} \times \cdots \times B_{b_d}}{\text{sup}} \sum_{i_1,...,i_d=1}^m \left| \Psi_{1,i_1}(z_1) \cdots \Psi_{d,i_d}(z_d) \right|.
\end{align*}

Hence, using (\ref{majphiij1}), we have 
\begin{equation} \label{unifpsi}
\underset{m \geq 1, \eta_{i_1,...,i_d} = \pm 1}{\text{sup}} \left\| \sum_{i_1,...,i_d=1}^m \eta_{i_1,...,i_d} \Psi_{1,i_1}(T_1)  \cdots \Psi_{d,i_d}(T_d) \right\|  < \infty.
\end{equation}

Likewise, by (\ref{majtildephiij}), we also have
\begin{equation} \label{uniftildepsi}
\underset{m \geq 1, \eta_{i_1,...,i_d} = \pm 1}{\text{sup}} \left\| \sum_{i_1,...,i_d=1}^m \eta_{i_1,...,i_d} \tilde{\Psi}_{1,i_1}(T_1)  \cdots \tilde{\Psi}_{d,i_d}(T_d) \right\|  < \infty.
\end{equation}

We now prove the quadratic $H^\infty(B_{\gamma_1} \times \cdots \times B_{\gamma_d})$ functional calculus property. Let $(\varphi_{j})_{j \in \mathcal{I}}$ be a finite family of functions of $H^\infty_{0,1} (B_{\gamma_1} \times \cdots \times B_{\gamma_d})$. Considering $Z = l^2_{\mathcal{I}}$ with its canonical basis denoted by $(e_{j})_{j \in \mathcal{I}}$, we define
\begin{equation} \label{defh}
h = \sum_{j \in \mathcal{I}} \varphi_{j} \otimes e_{j}
\end{equation}
regarded as an element of $H^\infty(B_{\gamma_1} \times \cdots \times B_{\gamma_d} ; Z) $.

Let now $(a_{i_1,...,i_d})$ be the family of $Z$ provided by (\ref{FMdvarialbes}) for $h$ defined in (\ref{defh}). We decompose every $a_{i_1,...,i_d}$ on basis $(e_{j})$ as
\begin{equation*}
a_{i_1,...,i_d} = \sum_{j \in \mathcal{I}} c_{(i_1,...,i_d);j} e_{j}, \qquad (i_1,...,i_d) \in (\mathbb{N}^*)^d.
\end{equation*}

For convenience, we will use notations $(i) = (i_1,...,i_d)$. Then we can write for every $j$ in $\mathcal{I}$ and $(\zeta_1,...,\zeta_d) \in B_{b_1} \times \cdots \times B_{b_d}$
\begin{equation} \label{FMvarphi}
\varphi_{j}(\zeta_1,...,\zeta_d) = \sum_{(i) \in \mathbb{N}^d} c_{(i);j} \Psi_{1,i_1}(\zeta_1)\tilde{\Psi}_{1,i_1}(\zeta_1) \cdots \Psi_{d,i_d}(\zeta_d)\tilde{\Psi}_{d,i_d}(\zeta_d).
\end{equation}

We now use that
\begin{equation*}
\left\| h \right\|_{\infty,B_{\gamma_1} \times \cdots \times B_{\gamma_d}}  = \left\| \left( \sum_{j \in \mathcal{I}} \left|\varphi_{j} \right|^2 \right)^\frac{1}{2} \right\|_{\infty,B_{\gamma_1} \times \cdots \times B_{\gamma_d}}
\end{equation*}
together with estimation (\ref{majalphaij1}) to say that
\begin{equation} \label{majaij}
\underset{(i_1,...,i_d) \in (\mathbb{N}^*)^d}{ \text{sup}} \left( \sum_{j \in \mathcal{I}} \left| c_{(i);j}  \right|^2 \right)^\frac{1}{2} \lesssim \left\| \left( \sum_{j \in \mathcal{I}} \left|\varphi_{j} \right|^2 \right)^\frac{1}{2} \right\|_{\infty,B_{\gamma_1} \times \cdots \times B_{\gamma_d}}.
\end{equation}

For any $j$ in $\mathcal{I}$ and any integer $m \geq 1$, we consider the function
\begin{equation*}
 h_{j; m} = \sum_{1 \leq i_1,...,i_d \leq m} c_{(i);j} \Psi_{1,i_1}\tilde{\Psi}_{1,i_1} \otimes \cdots  \otimes \Psi_{d,i_d}\tilde{\Psi}_{d,i_d},
 \end{equation*}
which belongs to $H^\infty_0(B_{\gamma_1} \times \cdots \times B_{\gamma_d}) $ and pointwise converges to $\varphi_{j}$ when $m \to \infty$, according to (\ref{FMvarphi}).

Fix now $(r_j)_{j \in \mathcal{I}}$ and $(r_{i_1,...,i_d})_{i_1,...,i_d \geq 1}$ two families of Rademacher independent variables on some probability space $(\Omega_0,\mathbb{P})$. Let $x$ in $X$. By the Khintchine-Kahane inequality, we have

\begin{align} \label{inegKK}
& \left\| \sum_{j \in \mathcal{I}} r_{j}  \otimes h_{j;m}(T_1,...,T_d) x \right\|_{\text{Rad}(X)} \\
\nonumber & = \left( \int_{\Omega_0} \left\| \sum_{(i),j} c_{(i);j} r_{j}(\omega) \Psi_{1,i_1}(T_1) \tilde{\Psi}_{1,i_1}(T_d) \cdots \Psi_{d,i_d}(T_d) \tilde{\Psi}_{d,i_d}(T_d) x \right\|^2 d\mathbb{P}(\omega) \right)^\frac{1}{2}\\
\nonumber & \lesssim \int_{\Omega_0} \left\| \sum_{(i)} \Psi_{1,i_1}(T_1) \cdots \Psi_{d,i_d}(T_d) \left( \sum_{j} c_{(i);(j)} r_{(j)}(\omega)  \tilde{\Psi}_{1,i_1}(T_d) \cdots  \tilde{\Psi}_{d,i_d}(T_d) x \right) \right\| d\mathbb{P}(\omega),
\end{align}
where indexes are such that $1 \leq i_1,...,i_d \leq m$ and $j$ in $\mathcal{I}$.

For any $(x_{i_1,...,i_d})$, we note that
\begin{align*}
& \sum_{1 \leq i_1,...,i_d \leq m} \Psi_{1,i_1}(T_1) \cdots \Psi_{d,i_d}(T_d) x_{i_1,...,i_d} \\
& = \int_{\Omega_0} \left(\sum_{1 \leq i_1,...,i_d \leq m} r_{i_1,...,i_d}(\omega) \Psi_{1,i_1}(T_1) \cdots \Psi_{d,i_d}(T_d) \right) \left( \sum_{1 \leq i_1,...,i_d \leq m} r_{i_1,...,i_d}(\omega) x_{i_1,...,i_d} \right) d \mathbb{P}(\omega),
\end{align*}
using that the $r_{i_1,...,i_d}$'s are independent.

Hence, by (\ref{unifpsi}), we obtain
\begin{align*}
\left\| \sum_{1 \leq i_1,...,i_d \leq m} \Psi_{1,i_1}(T_1) \cdots \Psi_{d,i_d}(T_d) x_{i_1,...,i_d} \right\| \lesssim \int_{\Omega_0} \left\|\sum_{1 \leq i_1,...,i_d \leq m} r_{i_1,...,i_d}(\omega) x_{i_1,...,i_d} \right\| d \mathbb{P}(\omega).
\end{align*}

We apply it in (\ref{inegKK}) with 
\begin{equation*}
x_{i_1,...,i_d} = \sum_{j \in \mathcal{I}} c_{(i_1,...,i_d);j} r_{j}(\omega)  \tilde{\Psi}_{1,i_1}(T_d) \cdots  \tilde{\Psi}_{d,i_d}(T_d) x 
\end{equation*}
to have
\begin{align*}
& \left\| \sum_{j \in \mathcal{I}} r_{j}  \otimes h_{j;m}(T_1,...,T_d) x \right\|_{\text{Rad}(X)} \lesssim \left\| \sum_{(i);j} c_{(i);j} r_{j} \otimes r_{(i)}  \tilde{\Psi}_{1,i_1}(T_d) \cdots  \tilde{\Psi}_{d,i_d}(T_d) x \right\|_{\text{Rad}(\text{Rad}(X))}.
\end{align*}

We recall now that from \cite{KaW}, we have
\begin{align*}
\left\| \sum _{(i);j} z_{(i);j} r_{(i)} \otimes r_{j} \otimes y_{(i)} \right\|_{\text{Rad}(\text{Rad}(X))} \lesssim \underset{(i)}{\text{sup}} \left( \sum_{j} \left| z_{(i);j} \right|^2 \right)^\frac{1}{2} \left\| \sum _{(i)} r_{(i)} \otimes y_{(i)} \right\|_{\text{Rad}(X)}
\end{align*}
for any finite family of complex numbers $(z_{(i);j})$ and elements $(y_{(i)})$ of $X$.

We apply this last inequality with $y_{(i)} = y_{i_1,...,i_d} =  \tilde{\Psi}_{1,i_1}(T_d) \cdots  \tilde{\Psi}_{d,i_d}(T_d) x $ together with (\ref{uniftildepsi}), (\ref{majaij}) and (\ref{majtildephiij}) to obtain
\begin{equation} \label{estquadrahmRitt}
\left\| \sum_{j \in \mathcal{I}} r_{j}  \otimes h_{j;m}(T_1,...,T_d) x \right\|_{\text{Rad}(X)} \lesssim \left\| x \right\| \left\| \left( \sum_{j \in \mathcal{I}} \left|\varphi_{j} \right|^2 \right)^\frac{1}{2} \right\|_{\infty,B_{\gamma_1} \times \cdots \times B_{\gamma_d}}.
\end{equation}

The final step of the proof uses approximation arguments. We use the same argument as in the end of the proof of \cite[Proposition 3.2]{ArLM}. We reproduce it here for convenience.

The inequality (\ref{estquadrahmRitt}) holds true when $(T_1,\ldots,T_d)$ is replaced by $(rT_1,\ldots,rT_d)$ for any $r\in(0,1)$.
Further, we know that $(h_{j;m})_{m\geq 1}$ is a bounded sequence of the space $H^\infty_0(B_{\gamma_1} \times \cdots \times B_{\gamma_d})$. Moreover, the sequence $(h_{j;m})_{m\geq 1}$ converges pointwise to $\varphi_{j}$. Hence applying Lebesgue's dominated convergence Theorem twice we have for any $j$
$$
\lim_{m\to\infty} h_{j;m}(rT_1,\ldots,rT_d) = \varphi_{j}(rT_1,\ldots,rT_d)
$$
for any $r\in(0,1)$, and
$$
\lim_{r\to 1^{-}} \varphi_{j}(rT_1,\ldots,rT_d) = \varphi_{j}(T_1,\ldots,T_d).
$$
We therefore deduce from (\ref{estquadrahmRitt}) that
\begin{equation*} 
\left\| \sum_{j \in \mathcal{I}} r_{j}  \otimes \varphi_{j}(T_1,...,T_d) x \right\|_{\text{Rad}(X)} \lesssim \left\| x \right\| \left\| \left( \sum_{j \in \mathcal{I}} \left|\varphi_{j} \right|^2 \right)^\frac{1}{2} \right\|_{\infty,B_{\gamma_1} \times \cdots \times B_{\gamma_d}},
\end{equation*}
which concludes the proof. $\square$

The proof of Proposition \ref{Quadrasquare} will require the following lemma.

\begin{lemma} \label{LemRitt}
Let $\alpha>0$ and $\gamma \in (0,\frac{\pi}{2})$. There exists a constant $C>0$ such that for any $z$ in $B_\gamma$, we have
\begin{equation} \label{unifborne}
 \left| 1-z \right|^{2\alpha} \sum_{k=1}^{+\infty} k^{2\alpha-1} (|z|^2)^{k-1} \leq C.
\end{equation} 

\end{lemma} 

\textit{Proof of Lemma \ref{LemRitt} : } We fix $x$ in $[0,1)$ and $\alpha>0$. Let $f : t \mapsto t^{2\alpha-1}x^{t-1}$ from $(0,\infty)$ into $\mathbb{R}$. If $2\alpha - 1 >0$, $f$ is incresing on $(0,\frac{2\alpha-1}{-\log(x)})$ and $f$ is decreasing on $( \frac{2\alpha-1}{-\log(x)},\infty) $.

Using a comparison test, we obtain an estimate
\begin{equation} \label{compserieintegral}
\sum_{k=1}^{\infty} k^{2\alpha-1} x^{k-1}  \leq C \int_0^\infty t^{2\alpha - 1} x^{t-1} dt,
\end{equation}
where the constant $C>0$ does not depend on $x$.

Changing variable $u = -\log(x) t$ in (\ref{compserieintegral}), we obtain
\begin{equation} \label{estserieintegral}
\sum_{k=1}^{\infty} k^{2\alpha-1} x^{k-1}  \lesssim \dfrac{\Gamma(2\alpha)}{x(-\log(x))^{2\alpha}},
\end{equation}
where $\Gamma(y) = \int_0^\infty u^{y-1} e^{-u} du, y >0$.

If we now take $z$ in a Stolz domain $B_\gamma$, with $z$ in a neighbourhood of $1$, we can apply (\ref{estserieintegral}) to have
\begin{equation*}
\left| 1-z \right|^{2\alpha} \sum_{k=1}^{n} k^{2\alpha-1} (|z|^2)^{k-1} \lesssim \dfrac{\Gamma(2\alpha)\left| 1-z \right|^{2\alpha}}{|z|^2 \left(-\log\left(\left|z\right|^2\right)\right)^{2\alpha}} \lesssim \dfrac{\Gamma(2\alpha)\left| 1-z \right|^{2\alpha}}{ 2^{2\alpha} \left|z\right|^2 (1-\left| z\right|)^{2\alpha}} \lesssim \left( \dfrac{\left|1-z\right|}{1-\left| z \right|}\right)^{2\alpha}.
\end{equation*}

We know that there exists a constant $A>0$ such that for every $\omega \in B_\gamma \setminus \left\lbrace 1 \right\rbrace$, we have $\frac{\left|1-\omega\right|}{1-\left| \omega \right|} \leq A$. Thus, $ \left| 1-z \right|^{2\alpha} \sum_{k=1}^{n} k^{2\alpha-1} (|z|^2)^{k-1} $ is uniformly bounded for $z \in B_\gamma$ and $n \in \mathbb{N}^*$.

If $2\alpha - 1<0$, the function $f : t \mapsto t^{2\alpha-1}x^{t-1}$ is decreasing on $\mathbb{R}^{+*}$ and one can proceed in the same way to obtain (\ref{compserieintegral}) and (\ref{estserieintegral}). $\square$

\textit{Proof of Proposition \ref{Quadrasquare} :} Let $T=(T_1,...,T_d)$ be a $d$-tuple of commuting Ritt operators and $\alpha = (\alpha_1,...,\alpha_d)$ in $(\mathbb{R}_+^*)^d$. In all this proof, we will use notation 
\begin{equation*} \label{xk}
x_{k_1,...,k_d} = \prod_{j=1}^d k_j^{\alpha_j-\frac{1}{2}} T_j^{k_j-1}(I_X-T_j)^{\alpha_j}(x),
\end{equation*}
in order to rewrite (\ref{fonccarreRitt}) as
\begin{equation}
\norm{x}_{T,\alpha} = \left\| \sum_{k_1,...,k_d \geq 1} r_{k_1,...,k_d} \otimes x_{k_1,...,k_d}\right\|_{\text{Rad}((\mathbb{N}^*)^d;X)}.
\end{equation}

As $X$ does not contain $c_0$, Theorem \ref{Kwapien} ensures that we only need to prove the existence of a constant $K>0$ such that for any integer $n \geq 1$ we have
\begin{equation} \label{partialsum}
\left\| \sum_{1 \leq k_1,...,k_d \leq n} r_{k_1,...,k_d} \otimes  x_{k_1,...,k_d} \right\|_{\text{Rad}((\mathbb{N}^*);X)} \leq K \left\| x \right\|.
\end{equation}

To have inequality (\ref{partialsum}), we apply (\ref{estquadraRitt}) with the functions 
\begin{equation*}
\varphi_{k_1,...,k_d}(z_1,...,z_d) = k_1^{\alpha_1 - \frac{1}{2}} \cdots k_d^{\alpha_d - \frac{1}{2}} z_1^{k_1-1}(1-z_1)^{\alpha_1} \cdots z_d^{k_d-1}(1-z_d)^{\alpha_d}.
\end{equation*}

It is clear that all these functions belong to $H^\infty_{0,1}(B_{\gamma_1} \times \cdots \times B_{\gamma_d})$ for any $\gamma_1,...,\gamma_d$ in $(0,\frac{\pi}{2})$. Moreover, definition of fractional power calculus says that
\begin{equation*}
\varphi_{k_1,...,k_d}(T_1,...,T_d) = k_1^{\alpha_1 - \frac{1}{2}} \cdots k_d^{\alpha_d - \frac{1}{2}} T_1^{k_1-1}(I_X-T_1)^{\alpha_1} \cdots z_d^{k_d-1}(I_X-T_d)^{\alpha_d}.
\end{equation*}

As $(T_1,...,T_d)$ has a quadratic $H^\infty(B_{\gamma_1} \times \cdots \times B_{\gamma_d})$ functional calculus by Proposition \ref{Calcquadra} we have an estimate 
\begin{align*}
\left\| \sum_{1 \leq k_1,...,k_d \leq n} r_{k_1,...,k_d} \otimes  \varphi_{k_1,...,k_d}(T_1,...,T_d) (x) \right\|_{\text{Rad}(X)} \lesssim \\ 
 \left\| x \right\| \left\| \left( \sum_{1 \leq k_1,...,k_d \leq n} \left| \varphi_{k_1,...,k_d} \right|^2 \right)^\frac{1}{2} \right\|_{\infty,B_{\gamma_1} \times \cdots \times B_{\gamma_d}}.
\end{align*}
Next we study the right hand side. Let $(z_1,...z_d) \in B_{\gamma_1} \times \cdots \times B_{\gamma_d}$. We write
\begin{align*} \label{regr}
\sum_{1 \leq k_1,...,k_d \leq n} \left| \varphi_{k_1,...,k_d}(z_1,...,z_d) \right|^2 = \prod_{j=1}^d \left[ \left| 1 - z_j \right|^{2\alpha_j}  \left( \sum_{k_j=1}^n k_j^{2\alpha_j-1} (|z_j|^{2})^{k_j-1} \right) \right],
\end{align*}

According to Lemma \ref{LemRitt}, there exists constants $C_1,...,C_d >0$ such that for any $j$ in $\left\lbrace 1,...,d \right\rbrace $, $z_j$ in $B_{\gamma_j}$ and $n$ in $\mathbb{N}^*$, we have
\begin{equation*}
\left| 1 - z_j \right|^{2\alpha_j}  \left( \sum_{k_j=1}^n k_j^{2\alpha_j-1} (|z_j|^{2})^{k_j-1} \right) \leq C_j.
\end{equation*}

Combining these inequalities, we see that $\prod_{j=1}^d \left[ \left| 1 - z_j \right|^{2\alpha_j}  \left( \sum_{k_j=1}^n k_j^{2\alpha_j-1} (|z_j|^{2})^{k_j-1} \right) \right]$ is uniformly bounded on $z_j \in B_{\gamma_j}$, $j=1,...,d$ and $n \in \mathbb{N}^*$, that is there exists a constant $K>0$ such that 
\begin{equation*}
\left\| \left( \sum_{1 \leq k_1,...,k_d \leq n} \left| \varphi_{k_1,...,k_d} \right|^2 \right)^\frac{1}{2} \right\|_{\infty,B_{\gamma_1} \times \cdots \times B_{\gamma_d}} \leq K,
\end{equation*}
with $K$ not depending on $n$. We finally obtain (\ref{partialsum}) and this suffices to prove the result. $\square$

\section{From square functions to $H^\infty$ joint functional calculus} \label{SquareJoint}

In this section, we show that square function estimates imply the $H^\infty$ joint functional calculus property for $R$-Ritt commuting operators. Let us recall this notion.

We take $(r_k)_{k \geq 1}$ a sequence of independent Rademacher variables. Let $X$ be a Banach space and let $E \subset B(X)$ be a set of bounded operators on $X$. We say that $E$ is $R$-bounded if there exists a constant $C \geq 0$ such that for any finite family $(T_k)_{1 \leq k \leq n}$ of $E$, $n \in \mathbb{N}^*$ and any finite family $(x_k)_{1 \leq k \leq n}$ of $X$ we have
\begin{equation} \label{Rborne}
\left\| \sum_{k=1}^n r_k \otimes T_k(x_k) \right\|_{\text{Rad}(X)} \leq C \left\| \sum_{k=1}^n r_k \otimes x_k \right\|_{\text{Rad}(X)}. 
\end{equation}

In this case, we let $\mathcal{R}(E)$ be the infimum of all the $C$ verifying (\ref{Rborne}).
We will also use the notion of $\gamma$-boundedness. If  $(g_k)_{k \geq 1}$ is a sequence of complex-valued independent standard Gaussian variables on some probability space $\Omega_0$, we say that $E$ is $\gamma$-bounded if there exists a constant $C \geq 0$ such that for any finite family $(T_k)_{1 \leq k \leq n}$ of $E$, $n \in \mathbb{N}^*$ and any finite family $(x_k)_{1 \leq k \leq n}$ of $X$ we have
\begin{equation*} 
\left\| \sum_{k=1}^n \gamma_k \otimes T_k(x_k) \right\|_{L_2(\Omega_0;X)} \leq C \left\| \sum_{k=1}^n \gamma_k \otimes x_k \right\|_{L_2(\Omega_0;X)}. 
\end{equation*}

Every $\gamma$-bounded subset of $B(X)$ is $R$-bounded. In the case where $X$ has finite cotype, a subset of $B(X)$ is $\gamma$-bounded if and only if it is $R$-bounded.

Here is a simple fact on products of $R$-bounded families (see e.g \cite[Proposition 8.1.19]{Hyt2}).

\begin{lemma} \label{lem23}
Let $X$ be a Banach space and $A_1,...,A_d$ be subsets of $B(X)$ such that every $A_k$ is $R$-bounded, $k=1,...,d$. Let 
\begin{equation*}
A_1 \cdots A_d = \left\lbrace S_1 \cdots S_d : S_k \in A_k, k=1,...,d \right\rbrace.
\end{equation*} 
Then $A_1 \cdots A_d$ is $R$-bounded with $\mathcal{R}(A_1 \cdots A_d) \leq \mathcal{R}(A_1) \cdots \mathcal{R}(A_d) $.
\end{lemma}

We will also use the following result, which is taken from \cite[8.5.2]{Hyt2} and is justified by the stability of $R$-boundedness property over closed convex combinations.

\begin{lemma} \label{lem22}
Let $X$ be a Banach space and $(\Omega,\mu)$ a measure space. Let $E$ be a $R$-bounded subset of $B(X)$ and $K>0$.  For every measurable function $F : \Omega \to E$, we let
\begin{equation*} \label{prodRborne}
E_K = \left\lbrace \int_{\Omega} h(\omega)F(\omega) d\mu(\omega) /  h \in L^1(\Omega ; \mu),   \int_{\Omega}\left| h(\omega)\right| d\mu(\omega) \leq K \right\rbrace.
\end{equation*}  

Then $E_K$ is $R$-bounded with $\mathcal{R}(E_K) \leq 2K \mathcal{R}(E) $.
\end{lemma}

A Ritt operator $T$ on $X$ is called $R$-Ritt provided that the two sets
\begin{equation*}
\left\lbrace T^n, n \in \mathbb{N} \right\rbrace, \qquad \left\lbrace n(T^n-T^{n-1}), n \geq 1 \right\rbrace
\end{equation*}
are $R$-bounded. In this case, there exists $\delta$ in $(0,\frac{\pi}{2})$ such that the set
\begin{equation*}
\left\lbrace (\lambda-1) R(\lambda,T) : z \in \mathbb{C} \setminus \overline{B_\delta} \right\rbrace
\end{equation*}
is $R$-bounded. We call $R$-type of $T$ the infimum of all $\delta$ for which this holds true (see \cite{AL} for details).

We now give a generalisation of the ergodic decomposition. Recall that for any power bounded operator $T$ acting on some reflexive Banach space $X$, we have $X = \text{Ker}(I_X-T) \oplus \overline{\text{Ran}(I_X-T)}$. Take now another power bounded operator $S$ acting on $X$ and commuting with $T$. As $S$ is a power bounded operator acting on subspaces $\text{Ker}(I_X-T)$ and $\overline{\text{Ran}(I_X-T)}$ which are reflexive, we can write ergodic decompositions on these spaces to have
\begin{equation*}
\text{Ker}(I_X-T) = (\text{Ker}(I_X-T) \cap \text{Ker}(I_X-S)) \oplus (\text{Ker}(I_X-T) \cap \overline{\text{Ran}(I_X-S)}),
\end{equation*}
\begin{equation*}
\overline{\text{Ran}(I_X-T)} = (\overline{\text{Ran}(I_X-T)} \cap \text{Ker}(I_X-S)) \oplus (\overline{\text{Ran}(I_X-T)} \cap \overline{\text{Ran}(I_X-S)})
\end{equation*}
and we obtain 
\begin{align*}
X =  &(\text{Ker}(I_X-T) \cap \text{Ker}(I_X-S)) \oplus (\text{Ker}(I_X-T) \cap \overline{\text{Ran}(I_X-S)}) \\
     & \oplus (\overline{\text{Ran}(I_X-T)} \cap \text{Ker}(I_X-S)) \oplus (\overline{\text{Ran}(I_X-T)} \cap \overline{\text{Ran}(I_X-S)}).
\end{align*}

Let now $T_1,...,T_d$ be commuting power bounded operators acting on reflexive Banach space $X$. Looking at the preceding discussion and using induction, one can see that we have the decomposition 
\begin{equation} \label{XLambda}
X = \bigoplus_{\Lambda \subset \left\lbrace 1,...,d \right\rbrace} X_\Lambda
\end{equation}
where we let
\begin{equation} \label{Ergodic}
X_\Lambda = \left[ \bigcap_{i \in \Lambda} \overline{\text{Ran}(I-T_i)} \right] \bigcap \left[ \bigcap_{i \notin \Lambda} \text{Ker}(I-T_i) \right]
\end{equation}
for any subset $\Lambda$ of $\left\lbrace 1,...,d \right\rbrace$.


\begin{theorem} \label{th21}

Let $X$ be a reflexive Banach space such that $X$ and $X^*$ have finite cotype. Let $T=(T_1,\ldots,T_d)$ be commuting Ritt operators on $X$ such that every $T_k$ is $R$-Ritt of $R$-type $\delta_k \in (0,\frac{\pi}{2})$ for $k=1,...,d$. Suppose that there exists a constant $C>0$ such that for any subset $\Lambda$ of $\{1,...,d\}$, there exist $\alpha_\Lambda = (\alpha_{k})_{k \in \Lambda}$ and $\beta_\Lambda = (\beta_{k})_{k \in \Lambda}$ in $(\mathbb{R}_+^*)^\Lambda$ such that
\begin{equation} \label{estcarre1}
\left\| x \right\|_{T,\alpha_\Lambda} \leq C \left\| x \right\|, \qquad x \in X,
\end{equation}

\begin{equation} \label{estcarre2}
\left\| y \right\|_{T^*,\beta_\Lambda} \leq C \left\| y \right\|, \qquad y \in X^*.
\end{equation}

Then $(T_1,...,T_d)$ admits a $H^\infty(B_{\gamma_1} \times \cdots \times B_{\gamma_d})$ joint functional calculus for every $\gamma_k \in (\delta_k , \frac{\pi}{2})$.

\end{theorem}  

\begin{remark}
We note that Theorem \ref{th21} and Theorem \ref{Thsquarefunc} are new results, even for a single operator $T$, when $\alpha \notin \mathbb{N}^*$.
\end{remark}

Firstly, we prove the equivalence between square functions on reflexive spaces.

\begin{theorem} \label{equivfonccarre}

Let $X$ be a reflexive Banach space with finite cotype. Let $T=(T_1,...,T_d)$ be a $d$-tuple of Ritt operators such that every $T_k$ is an $R$-Ritt operator, $k=1,...,d$. Let $\alpha = (\alpha_1,...,\alpha_d)$ and $\beta=(\beta_1,...\beta_d)$ be $d$-tuples of $(\mathbb{R}_+^*)^d$. Then $\left\| \cdot \right\|_{T,\alpha} $ and $\left\| \cdot \right\|_{T,\beta} $ are equivalent, that is there exists a constant $C >0$ such that
\begin{equation}
C^{-1} \left\| x \right\|_{T,\alpha} \leq \left\| x \right\|_{T,\beta} \leq C \left\| x \right\|_{T,\alpha}.
\end{equation}
\end{theorem}

\textit{Proof of Theorem \ref{equivfonccarre}}. This proof uses many ideas of the one of \cite[Theorem 3.3]{AL}. Let $\gamma_1,...,\gamma_d$ be positive numbers such that every $N_j = \alpha_j + \gamma_j$ is a positive integer, $j=1,...,d$. For every integers $k \geq 1$ and $j=1,...,d$, we let
\begin{equation*}
c_{k,j} = \frac{k(k+1) \cdots (k+N_j-2)}{k^{\alpha_j-\frac{1}{2}}}
\end{equation*}
if $N_j \geq 2$ and $c_{k,j} = \frac{1}{k^{\alpha_j-\frac{1}{2}}}$ if $N_j = 1$.
By \cite[Theorem 3.3]{AL}, we have for $j=1,...,d$ and $r$ in $(0,1)$
\begin{equation*} 
\sum_{k=1}^{\infty} c_{k,j} k^{\alpha_j-\frac{1}{2}} (rT_j)^{2k-2}(I_X - (rT_j)^2)^{N_j} = (N_j-1)! I_X
\end{equation*}
the series being absolutely convergent and as every $(I_X + r T_j)^{N_j}$ is invertible
\begin{equation} \label{identite1}
\sum_{k=1}^{\infty} c_{k,j} (rT_j)^{k-1}(I_X - rT_j)^{\gamma_j} k^{\alpha_j-\frac{1}{2}} (rT_j)^{k-1}(I_X - rT_j)^{\alpha_j} = (N_j-1)! (I_X + r T_j)^{-N_j}.
\end{equation}

Now let for any intergers $k_1,...,k_d,m_1,...,m_d \geq 1$ , $r$ in $(0,1)$ and $j=1,...,d$,
\begin{equation*}
S_{j,m_j,k_j} = m_j^{\beta_j-\frac{1}{2}} (rT_j)^{m_j+k_j-2} (I_X-rT_j)^{\beta_j+\gamma_j},
\end{equation*} 
\begin{equation*}
R_{j,k_j} =  k_j^{\alpha_j-\frac{1}{2}} (rT_j)^{k_j-1} (I_X-rT_j)^{\alpha_j}.
\end{equation*}

Then we define for any intergers $m_1,...,m_d \geq 1$, $r$ in $(0,1)$ and $x$ in $X$

\begin{equation*}
y_{(m_1,...,m_d)}(r) = \left[ \prod_{j=1}^d (N_j-1)! (I_X + r T_j)^{-N_j} m^{\beta_j - \frac{1}{2}} (rT_j)^{m_j-1} (I-rT_j)^{\beta_j} \right] x.
\end{equation*}

As every series in (\ref{identite1}) is absolutely convergent, one can see that
\begin{equation*}
y_{(m_1,...,m_d)}(r) = \sum_{k_1,...,k_d=1}^{\infty} c_{k_1,1} \cdots c_{k_d,d} S_{1,m_1,k_1} \cdots S_{d,m_d,k_d} R_{1,k_1} \cdots R_{d,k_d} x.
\end{equation*}

For any integer $n \geq 1$, define the partial sum
\begin{equation*}
y_{(m_1,...,m_d);n}(r) = \sum_{k_1,...,k_d=1}^{n} c_{k_1,1} \cdots c_{k_d,d} S_{1,m_1,k_1} \cdots S_{d,m_d,k_d} R_{1,k_1} \cdots R_{d,k_d} x.
\end{equation*}

Next we consider the square functions as follows. For any integers $k_1,...,k_d \geq 1$ we let
\begin{equation*}
x_{(k_1,...,k_d);(\alpha_1,...,\alpha_d)} = k_1^{\alpha_1-\frac{1}{2}} \cdots k_d^{\alpha_d-\frac{1}{2}} (rT_1)^{k_1-1} \cdots (rT_d)^{k_d-1} (I_X-rT_1)^{\alpha_1} \cdots (I_X-rT_d)^{\alpha_d} x
\end{equation*}
and similary for other integers $m_1,...,m_d \geq 1$
\begin{equation*}
x_{(m_1,...,m_d);(\beta_1,...,\beta_d)} = m_1^{\beta_1-\frac{1}{2}} \cdots m_d^{\beta_d-\frac{1}{2}} (rT_1)^{m_1-1} \cdots (rT_d)^{m_d-1} (I_X-rT_1)^{\beta_1} \cdots (I_X-rT_d)^{\beta_d} x.
\end{equation*}

By this way, square functions (\ref{fonccarreRitt}) may be written as
\begin{equation*}
\left\| x \right\|_{(rT_1,...,rT_d),(\alpha_1,...,\alpha_d)} = \left\| \sum_{k_1,...,k_d=1}^{\infty} r_{k_1,...,k_d} \otimes x_{(k_1,...,k_d);(\alpha_1,...,\alpha_d)} \right\|_{\text{Rad}(\mathbb{N}^d,X)}
\end{equation*}
and
\begin{equation} \label{SF1}
\left\| x \right\|_{(rT_1,...,rT_d),(\beta_1,...,\beta_d)} = \left\| \sum_{m_1,...,m_d=1}^{\infty} r_{m_1,...,m_d} \otimes x_{(m_1,...,m_d);(\beta_1,...,\beta_d)} \right\|_{\text{Rad}(\mathbb{N}^d,X)}.
\end{equation}

The aim is to have an estimate
\begin{equation} \label{inegsquarefunc}
\left\| x \right\|_{(rT_1,...,rT_d),(\beta_1,...,\beta_d)} \lesssim  \left\| x \right\|_{(rT_1,...,rT_d),(\alpha_1,...,\alpha_d)}.
\end{equation}

Let us study $  \left\| x \right\|_{(rT_1,...,rT_d),(\beta_1,...,\beta_d)}$. One can remark that
\begin{equation*}
y_{(m_1,...,m_d)}(r) = \left[ \prod_{j=1}^d (N_j-1)! (I_X + r T_j)^{-N_j} \right] x_{(m_1,...,m_d);(\beta_1,...,\beta_d)}.
\end{equation*}

As every set of operators $\left\lbrace (N_j-1)^{-1} (I_X + r T_j)^{N_j} : r \in (0,1) \right\rbrace$ is bounded for $j=1,...,d$, we can consider $\sum_{m_1,...,m_d=1}^{\infty} r_{m_1,...,m_d} \otimes y_{(m_1,...,m_d)(r)} $ in place of \\
$\sum_{m_1,...,m_d=1}^{\infty} r_{m_1,...,m_d} \otimes x_{(m_1,...,m_d);(\beta_1,...,\beta_d)} $ in (\ref{SF1}).

Fix now integers $n,M \geq 1$. Let $(m)=(m_1,...,m_d)$ and $(k) = (k_1,...,k_d)$. Considering that the summations run over all indexes $1 \leq m_1,...,m_d \leq M$ and $1 \leq k_1,...,k_d \leq n$, we have

\begin{align*}
\sum_{m_1,...,m_d} & r_{(m)} \otimes y_{(m);n}(r)\\
 								 & = \sum_{m_1,...,m_d} r_{(m)} \otimes \left( \sum_{k_1,...,k_d} c_{k_1,1} \cdots c_{k_d,d} S_{1,m_1,k_1} \cdots S_{d,m_d,k_d} R_{1,k_1} \cdots R_{d,k_d} x \right) \\
 								 & = \sum_{m_1,...,m_d} r_{(m)} \otimes \left( \sum_{k_1,...,k_d} \delta_{k_1,m_1} \cdots \delta_{k_d,m_d} W_{1,m_1,k_1} \cdots W_{d,m_d,k_d} R_{1,k_1} \cdots R_{d,k_d} x \right)
\end{align*} 
where we write $\delta_{k_j,m_j} = \frac{c_{k_j,j}m_j^{\beta_j-\frac{1}{2}}}{(m_j+k_j-1)^{\gamma_j+\beta_j}} $ in order to have $
c_{k,j} S_{j,m_j,k_j} = \delta_{k_j,m_j} W_{j,m_j,k_j}$, with 
\begin{equation*}
W_{j,m_j,k_j} = (m_j+k_j-1)^{\gamma_j+\beta_j} (rT_j)^{m_j+k_j-2} (I_X - rT_j)^{\beta_j+\gamma_j}, \qquad j=1,...,d.
\end{equation*}

All the sets $F_j = \left\lbrace W_{j,m_j,k_j} : m_j,k_j \geq 1, r \in (0,1]\right\rbrace$ are $R$-bounded according to \cite[Proposition 2.8]{AL} for $j=1,...,d$ and using Lemma \ref{lem23}, the set $F_1 \cdots F_d$ is also $R$-bounded (using notation (\ref{prodRborne})). Since $X$ has finite cotype, this last set is then $\gamma$-bounded. 
 
Moreover, the infinite matrices $[\delta_{k_j ,m_j}]_{k_j,m_j \geq 1}$ represent an element of $B(l^2)$ denoted by $h_j$, for $j=1,\ldots,d$. Then the operator $h_1 \otimes \cdots \otimes h_d$, whose matrix has coefficients $\delta_{k_1,m_1} \cdots \delta_{k_d,m_d}$ 
represents an element of $B \left( l^2 \overset{2}{\otimes} \cdots \overset{2}{\otimes} l^2 \right)$ and we have $\left\| h_1 \otimes \cdots \otimes h_d \right\| = \left\| h_1 \right\| \cdots \left\| h_d \right\|$. Moreover, the coefficients of $h_1 \otimes \cdots \otimes h_d$ belong to $[0,\infty)$.


With all these properties in hand and considering Gaussian averages in place of Rademacher averages, we can use \cite[Proposition 2.6]{AL} to say that we have an inequality of type
 
\begin{align*}
\left\| \sum_{1 \leq m_1,...,m_d \leq M}  r_{m_1,...,m_d} \otimes y_{(m_1,...,m_d);n}(r) \right\|_{\text{Rad}((\mathbb{N}^*)^d;X)} \\
\lesssim \left\| \sum_{k_1,...,k_d = 1}^{\infty}  r_{k_1,...,k_d} \otimes x_{(k_1,...,k_d);(\alpha_1,...,\alpha_d)} \right\|_{\text{Rad}((\mathbb{N}^*)^d;X)}
\end{align*}
which suffices to have the finiteness of square functions and to have inequality (\ref{inegsquarefunc}).

It remains to study the case where $r \to 1^-$ to obtain the result. Choose $\nu$ an integer such that $\nu \geq \alpha_j + 1$ and $\nu \geq \beta_j + 1$, $j=1,...,d$. Using \cite[Lemma 3.2 (3)]{AL}, the limit $r \to 1^-$ exists for $x$ in $\bigcap_{j=1}^d \text{Ran}((I-T_j)^\nu) $
and we have for such an element $x$
\begin{equation*}
\left\| x \right\|_{(T_1,...,T_d),(\beta_1,...,\beta_d)} \lesssim  \left\| x \right\|_{(T_1,...,T_d),(\alpha_1,...,\alpha_d)}.
\end{equation*}

Let now for some integers $p_1,...,p_d$ the operators 
\begin{equation*}
\theta_{p_k} = \frac{1}{p_k+1} \sum_{l_k=0}^{p_k} (I-T_k^{l_k}), \qquad k=1,...,d
\end{equation*}
and
\begin{equation*}
\Theta_{p_1,...,p_d} = \theta_{p_1} \cdots \theta_{p_d}.
\end{equation*}

It is clear that $ \Theta_{p_1,...,p_d}^\nu$ maps $X$ into $\bigcap_{j=1}^d \text{Ran}((I-T_j)^\nu) $. Thus, we have the uniform estimate
\begin{equation*}
\left\| \Theta_{p_1,...,p_d}^\nu(x) \right\|_{(T_1,...,T_d),(\beta_1,...,\beta_d)} \lesssim  \left\| \Theta_{p_1,...,p_d}^\nu(x) \right\|_{(T_1,...,T_d),(\alpha_1,...,\alpha_d)}
\end{equation*}
for any $x$ in $X$ and $p_1,...p_d \geq 0$. 

Since all the $T_k$ are power bounded, the sequences $(\theta_{p_k})_{p_k \geq 0}$ are bounded. Using \cite[Lemma 3.2 (1)]{AL}, we obtain
\begin{equation*}
\left\| \Theta_{p_1,...,p_d}^\nu(x) \right\|_{(T_1,...,T_d),(\beta_1,...,\beta_d)} \lesssim  \left\| x \right\|_{(T_1,...,T_d),(\alpha_1,...,\alpha_d)}
\end{equation*}
for any $x$ in $X$ and $p_1,...p_d \geq 0$. 

Further, we know that for any $x \in \bigcap_{j=1}^d \text{Ran}(I-T_j)$ we have $\Theta_{p_1,...,p_m}(x) \longrightarrow x$ and $\Theta_{p_1,...,p_m}^\nu (x) \longrightarrow x$ for $p_1,...,p_d \to \infty$. Considering finite sums with $q \geq 1$ in square functions with estimates 
\begin{align*} 
\left\| \sum_{ 1 \leq m_1,...,m_d \leq q} \prod_{j=1}^d m_j^{\beta_j-\frac{1}{2}} r_{m_1,...,m_d} \otimes \left( \prod_{j=1}^d T_j^{m_j} (I_X-T_j)^{\beta_j} \Theta_{p_1,...,p_m}^\nu(x)  \right)\right\|_{\text{Rad}((\mathbb{N}^*)^d;X)}  \\  \lesssim  \left\| x \right\|_{(T_1,...,T_d),(\alpha_1,...,\alpha_d)}
\end{align*}
and passing to the limit $p_1,...,p_d \to \infty$ in the left hand sum then letting $q \to \infty$ yields estimate $\left\| x \right\|_{(T_1,...,T_d),(\alpha_1,...,\alpha_d)} \lesssim  \left\| x \right\|_{(T_1,...,T_d),(\beta_1,...,\beta_d)}$ for every $x$ in $ \bigcap_{j=1}^d \text{Ran}(I-T_j)$.

To conclude the proof, take ergodic decomposition given by (\ref{XLambda}) and (\ref{Ergodic}).

We remark that square functions $\left\| x \right\|_{(T_1,...,T_d),(\alpha_1,...,\alpha_d)}$ and $\left\| x \right\|_{(T_1,...,T_d),(\beta_1,...,\beta_d)} $ vanish on every subspace $\text{Ker}(I-T_k)$, $k=1,...,d$, and then on every subspace $X_\Lambda$ where $\Lambda \neq \left\lbrace 1,...,d \right\rbrace $. This means that the estimate $\left\| x \right\|_{(T_1,...,T_d),(\beta_1,...,\beta_d)} \lesssim  \left\| x \right\|_{(T_1,...,T_d),(\alpha_1,...,\alpha_d)}$ on $ \bigcap_{j=1}^d \text{Ran}(I-T_j)$ suffice to have this one in all the space $X$, which ends the proof. $\square$

\textit{Proof of Theorem \ref{th21}}. Taking into account Theorem \ref{equivfonccarre}, we may and do assume that $ \alpha = \beta = (1,...,1)$. Let $\gamma_k \in (\delta_k,\frac{\pi}{2})$ for $k=1,...,d$. Let $x \in X$, $y \in X^*$ and $\varphi \in H_{0,1}(B_{\gamma_1} \times \cdots \times B_{\gamma_d})$.  The aim is to have an estimate
\begin{equation} \label{estHinfini}
\left| \langle \varphi(T_1,...,T_d) x , y \rangle \right| \lesssim \left\| x \right\| \left\| y \right\| \left\| \varphi \right\|_{\infty,B_{\gamma_1} \times \cdots \times B_{\gamma_d}}.
\end{equation} 

By \cite[Propostion 2.5]{ArLM}, it suffices to prove (\ref{estHinfini}) when $\varphi$ is a polynomial function on $d$ variables.

We suppose first that $\varphi$ is such a polynomial function of the form $\varphi(z_1,...,z_d) = (1-z_1) \cdots (1-z_d) \varphi_1(z_1,...,z_d)$ where $\varphi_1$ is another polynomial function.


Then we see that $\varphi(T_1,...,T_d)x \in \bigcap_{k=1}^d \text{Ran}(I-T_k)$ for every $x$ in $X$.

According to \cite[Lemma 7.2]{LM}, we have for every $k=1,..,d$ and $y$ in $\text{Ran}(I-T_k)$
\begin{equation*}
\sum_{i_k=1}^{\infty} i_k(i_k+1) T_k^{i_k-1}(I-T_k)^3 y = 2y.
\end{equation*}
We deduce that
\begin{align} \label{identite}
\sum_{i_d=1}^{\infty} \cdots \sum_{i_1=1}^{\infty} i_d(i_d+1) \cdots i_1(i_1+1) T_d^{i_d-1} \cdots T_1^{i_1-1} (I-T_d)^3 \cdots (I-T_1)^3 \varphi(T_1,...,T_d)x \\ \nonumber = 2^d \varphi(T_1,...,T_d) x, \qquad x \in X.
\end{align}

Consider the polynomial function $\psi$ defined by $\psi(z) = \frac{1}{2}(1+z+z^2)^3$ 
Also for any integers $i_1,...,i_d \geq 1$, we set
\begin{equation*}
f(i_1,...,i_d) = \prod_{k=1}^d (i_k+1) \prod_{k=1}^d T_k^{i_k-1}(I-T_k) \varphi(T_1,...,T_d),
\end{equation*}
\begin{equation*}
g(i_1,...,i_d) = \prod_{k=1}^d i_k^\frac{1}{2}\prod_{k=1}^d T_k^{i_k-1}(I-T_k), 
\end{equation*} 
\begin{equation*}
h(i_1,...,i_d) = \prod_{k=1}^d i_k^\frac{1}{2}\prod_{k=1}^d (T_k^*)^{i_k-1}(I-T_k^*)\psi(T_k^*).
\end{equation*}

For convenience, we will write only $\displaystyle{\sum_{i_1,...,i_d}}$ in place of $\displaystyle{\sum_{i_d=1}^{\infty} \cdots \sum_{i_1=1}^{\infty} }$, keeping the order of summation. Then it follows from (\ref{identite}) that for any $x$ in $X$ and $y$ in $X^*$ we have
\begin{equation} \label{dualite}
\langle \varphi(T_1,...,T_d) x , y \rangle  = \sum_{i_1,...,i_d} \langle f(i_1,...,i_d)g(i_1,...,i_d) x, h(i_1,...,i_d) y\rangle.
\end{equation}

Let us now consider independent Rademacher variables $(r_{i_1,...,i_d})_{(i_1,...,i_d) \in (\mathbb{N}^*)^d}$ and families $(x_{i_1,...,i_d})_{(i_1,...,i_d) \in (\mathbb{N}^*)^d}$ and $(y_{i_1,...,i_d})_{(i_1,...,i_d) \in (\mathbb{N}^*)^d}$ of $X$ and $X^*$ respectively. For any integers $N_1,...,N_d \geq 1$, writing $\displaystyle{\sum^{N_1,...,N_d}}$ in place of $\displaystyle{\sum_{i_1=1}^{N_1} \cdots \sum_{i=1}^{N_d}} $, the independance of the $r_{i_1,...,i_d}$ yields

\begin{equation} \label{independance}
\sum^{N_1,...,N_d} \langle x_{i_1,...,i_d} , y_{i_1,...,i_d} \rangle = \int_\Omega \left\langle \sum^{N_1,...,N_d} r_{i_1,...,i_d}(u)  x_{i_1,...,i_d}~~,~~ \sum^{N_1,...,N_d} r_{i_1,...,i_d}(u) y_{i_1,...,i_d} \right\rangle 
 d \mathbb{P}(u).
\end{equation}

We now let 
\begin{equation*}
S_{N_1,...,N_d} = \sum^{N_1,...,N_d} \langle f(i_1,...,i_d)g(i_1,...,i_d) x, h(i_1,...,i_d) y \rangle,
\end{equation*}
the partial sums of (\ref{dualite}) for any integers $N_1,...,N_d \geq 1$. Letting 
\begin{equation*}
x_{i_1,...,i_d} = f(i_1,...,i_d)g(i_1,...,i_d) x, \qquad y_{i_1,...,i_d} = h(i_1,...,i_d) y,
\end{equation*}
for any integers $i_1,...,i_d \geq 1$ and using the Cauchy-Schwarz inequality in (\ref{independance}), we obtain
\begin{align} \label{partsom}
\left| S_{N_1,...,N_d} \right| \leq & \left\| \sum^{N_1,...,N_d} r_{i_1,...,i_d} \otimes f(i_1,...,i_d) g(i_1,...,i_d) x \right\|_{\text{Rad}(X)} \times \\
\nonumber & \left\| \sum^{N_1,...,N_d} r_{i_1,...,i_d} \otimes h(i_1,...,i_d) y \right\|_{\text{Rad}(X^*)}.
\end{align}

We now prove that the family of operators $(f(i_1,...,i_d))$ is $R$-bounded. First, as $\varphi$ is a function of $H^\infty_0(B_{\gamma_1} \times \cdots \times B_{\gamma_d}) $, operators $f(i_1,...,i_d)$ are defined by
\begin{equation*} \label{fi}
f(i_1,...,i_d) = \left(\dfrac{1}{2\pi \text{i}} \right)^d \int_{\partial B_{(\gamma)}} \prod_{k=1}^d (i_k+1) \varphi(\lambda_1,...,\lambda_d) \prod_{k=1}^d \lambda_k^{i_k}(\lambda_k-1) \prod_{k=1}^{d} R(\lambda_k,T_k) d(\lambda),
\end{equation*}
where $\partial B_{(\gamma)} = \partial B_{\gamma_1} \times \cdots \times \partial B_{\gamma_d} $ and $d(\lambda) = d\lambda_1 \cdots d\lambda_d$.

We consider the subset of $B(X)$ 
\begin{equation} \label{ensRborne}
E = \left\lbrace (\lambda_1 -1) R(\lambda_1,T_1) \cdots (\lambda_1 -1) R(\lambda_1,T_1) : (\lambda_1,...,\lambda_d) \in \partial B_{\gamma_1} \times \cdots \times \partial B_{\gamma_d}\right\rbrace.
\end{equation}

Every $E_k = \left\lbrace (\lambda_k -1) R(\lambda_k,T_k) : \lambda_k \in \partial B_{\gamma_k} \right\rbrace$ is $R$-bounded as $T_k$ is $R$-Ritt of $R$-type $\delta_k$ and $\delta_k < \gamma_k$. Using Lemma \ref{lem23}, we know that $E$ is $R$-bounded.

Let us now use Lemma \ref{lem22}. Consider the function $F : (\lambda_1,...,\lambda_d) \mapsto (\lambda_1 -1) R(\lambda_1,T_1) \cdots (\lambda_1 -1) R(\lambda_1,T_1)$ from $\partial B_{(\gamma)} $ to $E$. We prove that $\chi_{i_1,...,i_d} : (\lambda_1,...,\lambda_d) \mapsto \prod_{k=1}^d (i_k+1) \lambda^{i_k} \varphi(\lambda_1,...,\lambda_d)$ is uniformly bounded in $L^1(\partial B_{(\gamma)}, |d(\lambda)|)$ when $(i_1,...,i_d) \in (\mathbb{N}^*)^d$. Recall that for any angle $\theta \in (0,\frac{\pi}{2})$, one can check that 
\begin{equation*}
\text{sup} \left\lbrace \int_{\partial B_\theta} (m+1) |z|^{m} |dz| , m \in \mathbb{N} \right\rbrace < \infty.
\end{equation*}
$\chi_{i_1,...,i_d}$ is therefore uniformly bounded in $L^1(\partial B_{(\gamma)}, |d(\lambda)|)$ with $\left\| h_{i_1,...,i_d} \right\|_{L^1(\partial B_{(\gamma)}, |d(\lambda)|)} \lesssim \left\| \varphi \right\|_{\infty, \partial B_{(\gamma)}} $. We apply Lemma \ref{lem22} to obtain that the family $(f(i_1,...,i_d))$ is $R$-bounded with 
\begin{equation} \label{constR}
\mathcal{R} \left\lbrace f(i_1,...,i_d) : (i_1,...,i_d) \in (\mathbb{N}^*)^d \right\rbrace \lesssim \left\| \varphi \right\|_{\infty, \partial B_{(\gamma)}}.
\end{equation}

We now use (\ref{constR}) together with (\ref{partsom}) to have
\begin{align*}
\left| S_{n_1,...,N_d} \right| \lesssim & \left\| \varphi \right\|_{\infty, \partial B_{(\gamma)}} \left\| \sum^{N_1,...,N_d} r_{i_1,...,i_d} \otimes g(i_1,...,i_d) x \right\|_{\text{Rad}(X)} \times \\
& \left\| \sum^{N_1,...,N_d} r_{i_1,...,i_d} \otimes h(i_1,...,i_d) y \right\|_{\text{Rad}(X^*)}. 
\end{align*}
Finally, we let successively $N_1 \to \infty$, ..., $N_d \to \infty $ to obtain 
\begin{equation*}
\left| \langle \varphi(T_1,...,T_d) x , y \rangle \right| \lesssim \left\| \varphi \right\|_{\infty, \partial B_{(\gamma)}} \left\| x \right\|_{(T_1,...,T_d),(1,...,1)} \left\| y \right\|_{(T_1^*,...,T_d^*),(1,...,1)}.
\end{equation*}

Then it suffices to use (\ref{estcarre1}) and (\ref{estcarre2}) with $\alpha=\beta=(1,...,1)$ to deduce that 
\begin{equation*}
\left| \langle \varphi(T_1,...,T_d) x , y \rangle \right| \lesssim \left\| \varphi \right\|_{\infty, B_{(\gamma)}} \left\| x \right\| \left\| y \right\| 
\end{equation*}
which provides $\left\| \varphi(T_1,...,T_d) \right\| \lesssim \left\| \varphi \right\|_{\infty, B_{(\gamma)}}$. 

To conclude the proof, one sees that we have the estimation (\ref{estHinfini}) for any function polynomial function $\varphi$ of type $\varphi(z_1,...,z_d) = \prod_{i \in \Lambda} (1-z_i) \varphi_1(z_1,...,z_d) $ where $\varphi_1$ is another polynomial function depending on variables $(z_i)_{i \in \Lambda}$ with $\Lambda$ a subset of $\left\lbrace 1,...,d \right\rbrace $. Indeed, arguments above holds verbatim using square functions $\left\| \cdot \right\|_{T,\alpha_\Lambda} $ and $\left\| \cdot \right\|_{T^*,\beta_\Lambda} $. Hence,  (\ref{estHinfini}) is verified for any polynomial function.
$\square$

\begin{remark}
Looking at the preceding proof, one sees that we do not need to suppose that $X$ is a reflexive space or $X$ has finite cotype in Theorem \ref{th21} if we only consider square functions estimates $\left\| x \right\|_{T,(1,...,1)} \lesssim  \left\| x \right\|$ and $\left\| y \right\|_{T^*,(1,...,1)} \lesssim  \left\| y \right\|$ for $x \in X$ and $y \in X^*$.
\end{remark}

Recall that for $p \in (1,\infty)$, $p \neq 2$, the noncommutative $L_p$-spaces do not have property $(\alpha)$ (see e.g \cite[Section 3]{ArLM} for a definition of this property). Thus, the results of \cite[Section 3]{ArLM} do not apply to this class of Banach space and in particular, noncommutative $L_p$-spaces do not have the joint functionnal calculus property (see \cite{LLLM}). The following result, which generalises \cite[Corollary 7.5]{LM}, gives a characterisation of the joint functional calculus of a $d$-tuple of Ritt operators on spaces having property $(\Delta)$ (see \cite{KW1} for the definition). This result applies to the noncommutative $L_p$-spaces. Note that we must appeal to Theorem \ref{Thsquarefunc} and Remark \ref{Remarque} to obtain the next corollary.

\begin{corollary}
Let $X$ be a Banach space with property $(\Delta)$. Let $T_1,...,T_d$ be commuting Ritt operators on $X$. The following two assertions are equivalent.
\begin{itemize}
\item[i)] $T=(T_1,...,T_d)$ admits a $H^\infty(B_{\gamma_1} \times \cdots \times B_{\gamma_d})$ joint functional calculus for some $\gamma_k \in (0,\frac{\pi}{2})$, $k=1,...,d$.
\item[ii)] Every $T_k$ is $R$-Ritt and  for every $\Lambda$ subset of $\left\lbrace 1,...,d \right\rbrace$ and for any $\alpha_\Lambda$ and $\beta_\Lambda$ in $ (\mathbb{R}_+^*)^\Lambda$, there exists a constant $C>0$ such that we have
\begin{equation*}
\left\| x \right\|_{T,\alpha_\Lambda} \leq C \left\| x \right\|, \qquad x \in X,
\end{equation*}
\begin{equation*}
\left\| y \right\|_{T^*,\beta_\Lambda} \leq C \left\| y \right\|, \qquad y \in X^*.
\end{equation*}
\end{itemize}

\end{corollary}

\section{From $H^\infty$ joint functional calculus to dilation} \label{Dilations}

In this section, we give an application of square functions in terms of dilation of a $d$-tuple of Ritt operators. The framework of this part is $K$-convex spaces. Let us recall some background on these spaces.

Let $(r_n)_{n \geqslant 1}$ be a sequence of independent Rademacher variables on any probability space $\Omega_0$. Denote by $R$ the orthogonal projection from $L_2(\Omega_0)$ onto $\text{Rad}(\mathbb{N}^*;\mathbb{C})$, the closed subspace spanned by all the $r_n$. 

A Banach space $X$ is called $K$-convex if the operator $R \otimes I_X$ defined a priori on the space $L_2(\Omega_0) \otimes X$ extends to a bounded operator on $L_2(\Omega_0;X)$ (see \cite{P5} and \cite[Section 6]{Ma}). Using Khintchine Kahane's inequalities, one can see that if $X$ is a $K$-convex space, the space $\text{Rad}_p(X)$ is complemented in $L_p(\Omega_0;X)$ for any $1<p<\infty$, namely $R \otimes I_X$ extends to a bounded operator on $L_p(\Omega_0;X)$. We call this extension ''canonical projection'' from $L_p(\Omega_0;X)$ onto $\text{Rad}_p(X)$.

Let now $Y$ be another Banach space and $U_1,...,U_d$ be commuting isomorphisms on $Y$. We say that $U=(U_1,...,U_d)$ admits a $C(\mathbb{T}^d)$ bounded functional calculus if there exists a constant $C \geq 1$ such that for any trigonometric polynomial of $d$ variables 
\begin{equation*}
\phi(s_1,...,s_d) = \sum_{n_1,...,n_d \in \mathbb{Z}^d} a_{n_1,...,n_d} e^{\text{i}(n_1 s_1 + \cdots + n_d s_d)},
\end{equation*} 
where $(a_{n_1,...,n_d})$ is a finite family of complex numbers, we have 
\begin{equation}
\left\| \phi(U_1,...,U_d) \right\| \leq C \text{sup} \left\lbrace \left| \phi(s_1,...,s_d) \right| : (s_1,...,s_d) \in \mathbb{R}^d \right\rbrace,
\end{equation}
where we let $\phi(U_1,...,U_d) = \sum a_{n_1,...,n_d} U_1^{n_1} \cdots U_d^{n_d}$.

As trigonometric polynomial functions of $d$ variables are dense in the space $C(\mathbb{T}^d)$ of all continuous functions on $\mathbb{T}^d$, $C(\mathbb{T}^d)$ bounded functional calculus property is equivalent to the existence of a unique unital bounded homomorphism $\omega : C(\mathbb{T}^d) \to B(X)$ such that $\omega(e(j)) = U_j$ where $e(j) : (s_1,...,s_d) \mapsto e^{\text{i}s_j}$.

\begin{theorem} \label{Dilation}
Let $X$ be a reflexive $K$-convex Banach space and $p$ in $(1,\infty)$. Let $T=(T_1,...,T_d)$ be a $d$-tuple of commuting Ritt operators on $X$. Suppose more that $T$ admits a $H^\infty(B_{\gamma_1} \times \cdots \times B_{\gamma_d})$ joint functional calculus for some $\gamma_1,...,\gamma_d$ in $(0,\frac{\pi}{2})$.

Then there exists a measure space $\Sigma$, a $d$-tuple of commuting isomorphisms $(U_1,...,U_d)$ on $L_p(\Sigma;X)$ admitting a $C(\mathbb{T}^d)$ bounded calculus and two bounded operators $J : X \to L_p(\Sigma;X)$ and $Q : L_p(\Sigma;X) \to X$ such that
\begin{equation} \label{dilationisom}
T_1^{n_1} \cdots T_d^{n_d} = Q U_1^{n_1} \cdots U_d^{n_d} J, \qquad (n_1,...,n_d) \in \mathbb{N}^d.
\end{equation}
\end{theorem}

\textit{Proof of Theorem \ref{Dilation} :} Throughout the proof, we let $\Lambda$ be a subset of $\left\lbrace 1,...,d \right\rbrace $, denoted by
\begin{equation} \label{ensemble}
\Lambda = \left\lbrace i_1,...,i_k \right\rbrace, \qquad i_1 < \cdots < i_k,
\end{equation}
where $|\Lambda| = k$. Let $(r_{(l_{i_1},...,l_{i_k})}) $ be a family of independent Rademacher variables indexed by $\mathbb{Z}^\Lambda$ for $\Lambda \neq \emptyset$ on some probability space $\Omega_\Lambda$. Consider then space $\text{Rad}_p(\mathbb{Z}^\Lambda;\mathbb{C})$, letting $\text{Rad}_p(\mathbb{Z}^\emptyset;\mathbb{C}) = \mathbb{C}$.

For any $j$ in $\left\lbrace 1,...,d \right\rbrace$ and $\Lambda \neq \emptyset$, we define the operator $v_{j,\Lambda}$ on the space $\text{Rad}_p(\mathbb{Z}^\Lambda;\mathbb{C})$ by $v_{j,\Lambda}(r_{(l_{i_1},...,l_{i_k})}) = r_{(l_{i_1},...,l_{i_k})} $ if $j \notin \Lambda$ and if $j=i_m \in \Lambda$, we let $v_{j,\Lambda}(r_{(l_{i_1},...,l_{i_m},...,l_{i_k})}) = r_{(l_{i_1},...,l_{i_m}-1,...,l_{i_k})} $. If $\Lambda = \emptyset$, let $v_{j,\emptyset} = I_\mathbb{C}$ for $j=1,...,d$.

Then $v_{j,\Lambda}$ is an isometric isomorphism of $\text{Rad}_p(\mathbb{Z}^\Lambda;\mathbb{C})$. Further, every $v_{j,\Lambda} \otimes I_X$ extends to a unique operator $V_{j,\Lambda}$ from $\text{Rad}_p(\mathbb{Z}^\Lambda;X) $ into itself, which is an isometric isomorphism too.

Next we define 
\begin{equation*}
V_j = \bigoplus_{\Lambda \subset \left\lbrace 1,...,d \right\rbrace} V_{j,\Lambda},
\end{equation*}
the direct sum of operators $V_{j,\Lambda}$,  from $\overset{p}{\bigoplus}_{\Lambda \subset \left\lbrace 1,...,d \right\rbrace} \text{Rad}_p(\mathbb{Z}^\Lambda;X)$ into $\overset{p}{\bigoplus}_{\Lambda \subset \left\lbrace 1,...,d \right\rbrace} \text{Rad}_p(\mathbb{Z}^\Lambda;X)$ defined by $V_j(\sum_{\Lambda \subset \left\lbrace 1,...,d \right\rbrace} x_\Lambda) = \sum_{\Lambda \subset \left\lbrace 1,...,d \right\rbrace} V_{j,\Lambda}(x_\Lambda)$.

Recall now (\ref{XLambda}) and similary $(X^*)_\Lambda = \left[ \cap_{i \in \Lambda} \overline{\text{Ran}(I-T_i^*)} \right] \bigcap \left[ \cap_{i \notin \Lambda} \text{Ker}(I-T_i^*) \right]$. As $X$ is reflexive, by (\ref{Ergodic}), we have 
\begin{equation*}
X = \bigoplus_{\Lambda \subset \left\lbrace 1,...,d \right\rbrace} X_\Lambda \qquad \text{and} \qquad X^* = \bigoplus_{\Lambda \subset \left\lbrace 1,...,d \right\rbrace} (X^*)_\Lambda.
\end{equation*}

For any $\Lambda \subset \left\lbrace 1,...,d \right\rbrace$, taking notation of (\ref{ensemble}), we know that $T_\Lambda=(T_{i_1},...,T_{i_k})$ has a joint functional calculus as a subfamily of $T$. Using Theorem \ref{Thsquarefunc} with the $|\Lambda|$-tuple $\alpha = (\frac{1}{2},...,\frac{1}{2})$, we can define a bounded operator $J_{\Lambda,1}$ using square functions. 
We let
\begin{equation*}
J_{\Lambda,1} : \begin{array}{ccl} X_\Lambda &\to & \text{Rad}_p(\mathbb{Z}^\Lambda;X) \\ x & \mapsto & \displaystyle{\sum_{l_{i_1},...,l_{i_k}\geq 1} r_{(l_{i_1},...,l_{i_k})} \otimes \prod_{s=1}^k T_{i_s}^{l_{i_s}-1} \left( I_X -T_{i_s} \right)^\frac{1}{2} x} \end{array}.
\end{equation*}

Since $H^\infty$ functional calculus passes to the adjoint, we have square function estimates for $(T_1^*,...,T_d^*)$ and we can define in the same way
\begin{equation*}
J_{\Lambda,2} : \begin{array}{ccl} (X^*)_\Lambda &\to & \text{Rad}_{p'}(\mathbb{Z}^\Lambda;X^*) \\ y & \mapsto & \displaystyle{\sum_{l_{i_1},...,l_{i_k}\geq 1} r_{(l_{i_1},...,l_{i_k})} \otimes \prod_{s=1}^k ((T_{i_s})^*)^{l_{i_s}-1} \left( I_X -(T_{i_s})^* \right)^\frac{1}{2} y} \end{array},
\end{equation*}
where $p'$ is the conjugate exponent of $p$ (verifying $\frac{1}{p} + \frac{1}{p'} = 1 $).

For any $x \in X_\Lambda$, $y \in (X^*)_\Lambda$ and $n_1,...,n_k$ integers, we compute 
\begin{align} \label{CalcJ}
\langle V_{i_1}^{n_1} \cdots V_{i_k}^{n_k} J_{\Lambda,1} x , J_{\Lambda,2} y \rangle &= \sum_{l_{i_1},...,l_{i_k} \geq 1} \left\langle \prod_{s=1}^k T_{i_s}^{l_{i_s}+n_s-1} \left( I_X -T_{i_s} \right)^\frac{1}{2} x,\prod_{s=1}^k ((T_{i_s})^*)^{l_{i_s}-1} \left( I_{X^*} -(T_{i_s})^* \right)^\frac{1}{2} y \right\rangle \\
 \nonumber & = \sum_{l_{i_1},...,l_{i_k} \geq 1} \left\langle \prod_{s=1}^k T_{i_s}^{n_s + 2(l_{i_s}-1)} \left( I_X -T_{i_s} \right) x, y \right\rangle \\
 \nonumber & = \sum_{l_{i_1},...,l_{i_k} \geq 1} \left\langle \prod_{s=1}^k T_{i_s}^{n_s + 2(l_{i_s}-1)} \left( I_X -T_{i_s}^2 \right) (I_X+T_{i_s})^{-1} x, y \right\rangle.
\end{align}

Now recall that for any $x$ in $\overline{\text{Ran}(I_X-T_i)}$, $i=1,...,d$ we have
\begin{equation*}
\sum_{k=1}^\infty T_i^{2(k-1)}(I_X - T_i^2)(x) = x.
\end{equation*}
Indeed, one has $\sum_{k=1}^P T_i^{2(k-1)}(I_X - T_i^2)(x) = x - T_i^{2P}(x)$ for any $x$ and $P \geq 1$. If $x = (I_X-T_i)z$, $z \in X$, we have $\left\| T_i^{2P}(I_X-T_i) z \right\| \lesssim \frac{\left\|z  \right\|}{2P} $ using Ritt condition and this term tends to $0$ as $P$ tends to $\infty$.

Hence, letting $S_\Lambda = \prod_{s=1}^k (I+T_{i_s})^{-1}$ and developping the last sum in (\ref{CalcJ}), we finally obtain
\begin{equation*}
\langle V_{i_1}^{n_1} \cdots V_{i_k}^{n_k} J_{\Lambda,1} x , J_{\Lambda,2} y \rangle = \langle T_{i_1}^{n_1} \cdots T_{i_k}^{n_k} S_\Lambda x , y \rangle.
\end{equation*}

Next, for any $j \notin \Lambda$, we see that $T_j(x) = x$,$(T_j)^*(y) = y$ for $x \in X_\Lambda$ and $y \in (X^*)_\Lambda$. Moreover, $V_j$ acts as identity operator on $\text{Rad}(\mathbb{Z}^\Lambda;X) $, so that for any integers $n_1,...,n_d$
\begin{equation} \label{duali}
\langle V_{1}^{n_1} \cdots V_{d}^{n_d} J_{\Lambda,1} x , J_{\Lambda,2} y \rangle = \langle T_{1}^{n_1} \cdots T_{d}^{n_d} S_\Lambda x , y \rangle, \qquad x \in X_\Lambda, \qquad y \in (X^*)_\Lambda.
\end{equation}

Note that for any $a$ in $\text{Ker}(I_X-T_i)$ and $b$ in $\overline{\text{Ran}(I_X-T_i^*)}$, we have clearly $\left\langle a,b \right\rangle = 0$, $i=1,...,d$. Thus, for any $a$ in $X_\Lambda$ and $b$ in $(X^*)_{\Lambda'}$ with $\Lambda\neq \Lambda'$, one has $\left\langle a,b \right\rangle = 0$. Since each $T_{i}$ maps $X_\Lambda$ into itself for any $i$ in $\left\lbrace 1,...,d \right\rbrace$ and $\Lambda \subset \left\lbrace 1,...,d \right\rbrace$, one obtains that for subsets $\Lambda \neq \Lambda'$, we have 
\begin{equation} \label{duali2}
\langle V_{1}^{n_1} \cdots V_{d}^{n_d} J_{\Lambda,1} x , J_{\Lambda',2} y \rangle = 0, \qquad x \in X_\Lambda, \qquad y \in (X^*)_{\Lambda'}.
\end{equation}

$X$ being a $K$-convex space, every space $\text{Rad}_p(\mathbb{Z}^\Lambda;X) $ is complemented in $L_p(\Omega_\Lambda;X)$. Thus, the space $W = \overset{p}{\oplus}_{\Lambda \subset \left\lbrace 1,...,d \right\rbrace} \text{Rad}_p(\mathbb{Z}^\Lambda;X)$ is complemented in $\overset{p}{\oplus}_{\Lambda \subset \left\lbrace 1,...,d \right\rbrace} L_p(\Omega_\Lambda;X)$, which is identified to $L_p(\Sigma;X)$ where $\Sigma = \sqcup_{\Lambda \subset \left\lbrace 1,...,d \right\rbrace} \Omega_\Lambda$ equipped with the sum measure. Then we write $L_p(\Sigma;X) = W \oplus E$. It is the same for 
$$W' = \overset{p'}{\oplus}_{\Lambda \subset \left\lbrace 1,...,d \right\rbrace} \text{Rad}_{p'}(\mathbb{Z}^\Lambda;X^*),$$
which is complemented in $L_{p'}(\Sigma;X^*)$.

We can define operators $J_0$ and $Q_0$ as follows. We let 
\begin{equation*}
J_0 = \oplus_{\Lambda \subset\left\lbrace 1,...,d \right\rbrace}  J_{\Lambda,1} (S_\Lambda)^{-1}
\end{equation*} 
from $X$ into $W$, then we let $J_2 = \oplus_{\Lambda \subset\left\lbrace 1,...,d \right\rbrace}  J_{\Lambda,2}$ from $X^*$ into $W'$ and finally $Q_0 = J_2^*$, which maps $(W')^*$ into $X$. Using that $W$ acts by duality on $W'$ and the natural complementation of $W$ in $L_p(\Sigma;X)$ as well as complementation of $W'$ in $L_{p'}(\Sigma;X^*)$ (which relies on the canonical projection), we identify $(W')^*$ to $W$, so that $Q_0$ maps $W$ into $X$.
Identities (\ref{duali}) and (\ref{duali2}) guarantee that we have 
\begin{equation} \label{DilRad}
T_1^{n_1} \cdots T_d^{n_d} = Q_0 V_1^{n_1} \cdots V_d^{n_d} J_0, \qquad n_1,...,n_d \in \mathbb{N}.
\end{equation}  

We let $J = \mathcal{J} \circ J_0$ where $\mathcal{J} : W \hookrightarrow L_p(\Sigma;X)$ is the inclusion, $Q = Q_0 \circ \Pi$ where $\Pi$ is the canonical projection from $L_p(\Sigma;X)$ onto $W$ and finally let $U_j = V_j \oplus I_E$, $j=1,\ldots,d$. Then every $U_j$ is an isomorphism of $L_p(\Sigma;X)$ for $j=1,...,d$ having with (\ref{DilRad})
\begin{equation*}
T_1^{n_1} \cdots T_d^{n_d} = Q U_1^{n_1} \cdots U_d^{n_d} J \qquad n_1,...,n_d \in \mathbb{N}.
\end{equation*} 

It is clear that $U_1,...,U_d$ are commuting operators and we obtain (\ref{dilationisom}).

It remains to prove that $U = (U_1,...,U_d)$ admits a $C(\mathbb{T}^d)$ bounded calculus. It suffices to prove this fact for $V = (V_1,...,V_d)$ and hence for any $(V_{1,\Lambda},...,V_{d,\Lambda})$ for any $\Lambda \subset \left\lbrace 1,...,d \right\rbrace$. Let $\phi = \sum a_{n_1,...,n_d} e^{in_1\cdot} \cdots e^{in_d \cdot}$ be a trigonometric polynomial function. Noting that $V_{j,\Lambda} = v_{j,\Lambda} \overline{\otimes} I_X$, we have $\phi(V_{1,\Lambda},...,V_{d,\Lambda}) = \phi(v_{1,\Lambda},...,v_{d,\Lambda}) \overline{\otimes} I_X$.

We take $(g_{l_{i_1},...,l_{i_k}}) $ a sequence of independant Gaussian variable on some probability space $\Xi$ and we can define the space $\text{G}_2(\mathbb{Z}^\Lambda)$ as the closed subspace of $L_2(\Xi)$ spanned by $(g_{l_{i_1},...,l_{i_k}}) $. Next we define $\tilde{v_{j,\Lambda}} $ on  $\text{G}_2(\mathbb{Z}^\Lambda)$ in the same manner as $v_{j,\Lambda}$ on $\text{Rad}_p(\mathbb{Z}^\Lambda) $. It is clear that $\tilde{v_{j,\Lambda}} $ is a unitary operator of the Hilbert space $\text{G}_2(\mathbb{Z}^\Lambda)$. According to the Spectral Theorem, we obtain 
\begin{equation} \label{SpTh}
\left\| \phi(\tilde{v_{1,\Lambda}},...,\tilde{v_{d,\Lambda}}) \right\|_{\text{G}_2(\mathbb{Z}^\Lambda)\to \text{G}_2(\mathbb{Z}^\Lambda)} = \text{sup} \left\lbrace \left| \phi(s_1,...,s_d) \right|, (s_1,...,s_d) \in \mathbb{R}^d \right\rbrace.
\end{equation}

Recall that for any bounded operator $\psi : \text{G}_2(\mathbb{Z}^\Lambda) \to \text{G}_2(\mathbb{Z}^\Lambda) $, the operator $\psi \otimes I_X: \text{G}_2(\mathbb{Z}^\Lambda) \otimes X \to \text{G}_2(\mathbb{Z}^\Lambda) \otimes X$ extends to a unique bounded operator $\psi \overline{\otimes} I_X : \text{G}_2(\mathbb{Z}^\Lambda;X) \to \text{G}_2(\mathbb{Z}^\Lambda;X)$ with $\left\| \psi \overline{\otimes} I_X \right\| = \left\| \psi \right\| $ (see \cite[Prop. 6.1.23]{Hyt2}).

As $X$ is a $K$-convex space, it has finite cotype. Thus, the Rademacher averages and Gaussian averages on $X$ are equivalent. More precisely, the spaces $\text{Rad}_p(\mathbb{Z}^\Lambda;X)$ and $\text{G}_2(\mathbb{Z}^\Lambda;X)$ are naturally isomorphic and we have therefore an equivalence
\begin{equation*} 
\left\| \phi(V_{1,\Lambda},...,V_{d,\Lambda}) \right\|_{\text{Rad}_p(\mathbb{Z}^\Lambda;X) \to \text{Rad}_p(\mathbb{Z}^\Lambda;X)} \simeq \left\| \phi(\tilde{v_{1,\Lambda}},...,\tilde{v_{d,\Lambda}}) \overline{\otimes} I_X \right\|_{\text{G}_2(\mathbb{Z}^\Lambda;X) \to \text{G}_2(\mathbb{Z}^\Lambda;X)}
\end{equation*}
and then
\begin{equation} \label{GaussRad}
\left\| \phi(V_{1,\Lambda},...,V_{d,\Lambda}) \right\|_{\text{Rad}_p(\mathbb{Z}^\Lambda;X) \to \text{Rad}_p(\mathbb{Z}^\Lambda;X)} \simeq \left\| \phi(\tilde{v_{1,\Lambda}},...,\tilde{v_{d,\Lambda}}) \right\|_{\text{G}_2(\mathbb{Z}^\Lambda) \to \text{G}_2(\mathbb{Z}^\Lambda)}.
\end{equation}

Combining  (\ref{SpTh}) and (\ref{GaussRad}), we obtain an inequality of type
\begin{equation*}
\left\| \phi(V_{1,\Lambda},...,V_{d,\Lambda}) \right\|_{\text{Rad}(\mathbb{Z}^\Lambda;X)} \lesssim \text{sup} \left\lbrace \left| \phi(s_1,...,s_d) \right|, (s_1,...,s_d) \in \mathbb{R}^d \right\rbrace,
\end{equation*}
which is exactly the property of $C(\mathbb{T}^d) $ bounded calculus. $\square$

\section{From dilation to $H^\infty$ joint functional calculus} \label{HRRitt}

The aim of this section is to have a converse property of Theorem \ref{Dilation}, which is stated as follows.

\begin{theorem} \label{Diljointcalc}

Let $X$ be a Banach space and $p \in (1,\infty)$. Let $(T_1,...,T_d)$ be a $d$-tuple of commuting operators acting on $X$ such that every $T_k$ is a $R$-Ritt operator, $k=1,...,d$. Suppose that there exist a measure space $\Sigma$, a $d$-tuple of commuting isomorphisms $(U_1,...,U_d)$ acting on $L_p(\Sigma;X)$ having a $C(\mathbb{T}^d)$ bounded calculus and two bounded operators $J : X \to L_p(\Sigma;X)$, $Q : L_p(\Sigma;X) \to X$ such that (\ref{dilationisom}) is verified. 

Then  there exist $b_1,...,b_d$ in $(0,\frac{\pi}{2})$ such that $(T_1,...,T_d)$ admits a $H^\infty(B_{b_1} \times \cdots \times B_{b_d})$ joint functional calculus.

\end{theorem}

To obtain the Theorem above, we generalise a result of Franks and McIntosh (\cite[Theorem 5.1]{FMI}). This result allows to reduce the domain for which a certain operator on an Hilbert space admits an $H^\infty$ functional calculus. We show that this admits a generalisation with Stolz domains for general Banach spaces provided that we consider $R$-Ritt operators.


\begin{theorem} \label{GeneFMI}

Let $X$ be a Banach space and $(T_1,...,T_d)$ a $d$-tuple of commuting operators acting on $X$ such that every $T_k$ is an $R$-Ritt operator, $k=1,...,d$. Suppose that $(T_1,...,T_d)$ is polynomially bounded, that is there exists a constant $C \geq 1$ such that for any polynomial function $h$ of $d$ variables we have
\begin{equation} \label{polbound}
\left\| h(T_1,...,T_d) \right\| \leq C \text{sup} \left\lbrace \left| h(z_1,...,z_d) \right| : (z_1,...,z_d) \in \mathbb{T}^d \right\rbrace.
\end{equation}

Then $(T_1,...,T_d)$ admits a $H^\infty(B_{b_1} \times \cdots \times B_{b_d})$ joint functional calculus for some $b_k$ in $(0,\frac{\pi}{2})$, $k=1,...,d$.
\end{theorem}

For any integer $d \geq 1$, we let $H^\infty_0(\mathbb{D}^d)$ be the algebra of all holomorphic functions $f$ on $\mathbb{D}^d$ such that there exist positive constants $c$ and $(s_i)_{1 \leq i \leq d}$ such that
\begin{equation*}
\left| f(z_1,...,z_d) \right| \leq c \prod_{i=1}^d \left|1-z_i\right|^{s_i}, \qquad (z_1,...,z_d) \in \mathbb{D}.
\end{equation*}

The proof of Theorem \ref{GeneFMI} above requires a technical result which relies on the Franks-McIntosh decomposition on Stolz domain (see \cite[Section 6]{ArLM}). 
The proposition below is a direct consequence of \cite[Proposition 3.4]{ArLM2}, changing variables $z$ into $1-z$. We omit the details.

\begin{proposition} \label{DecompFMI}
Let $\gamma \in (0,\frac{\pi}{2})$. There exist sequences $(\phi_i)_{i \geq 1}$, $(\varphi_i)_{i \geq 1}$, $(\theta_i)_{i \geq 1}$ and $(\psi_i)_{i \geq 1}$ of $H^\infty_0(\mathbb{D})$ such that

\begin{itemize}
\item[(i)] For every $i \geq 1$, we have $\phi_i = \theta_i \varphi_i \psi_i$;
\item[(ii)] There exists a constant $c >0$ such that for every $z$ in $\mathbb{D}$
\begin{equation*}
\sum_{i=1}^\infty \left| \varphi_i(z) \right| \leq c, \qquad \sum_{i=1}^\infty \left| \psi_i(z) \right| \leq c;
\end{equation*}
\item[(iii)] There exists a constant $e >0$ such that for every $i \geq 1$ 
\begin{equation*}
\int_{\partial B_\gamma} \frac{\left| \theta_i(z) \right|}{\left| 1 - z \right|} \left| dz \right| \leq e;
\end{equation*}
\item[(iv)] For every $z$ in $\mathbb{D}$, the series $\sum_{i \geq 1} \phi_i(z)$ absolutely converges and there exists a constant $c'$ such that
\begin{equation*}
\text{sup} \left\lbrace \sum_{i=1}^{\infty} \left| \phi_i(z) \right| : z \in \mathbb{D} \right\rbrace \leq c'.
\end{equation*}
Moreover, we have
\begin{equation*}
\sum_{i=1}^{\infty} \phi_i(z) = 1, \qquad z \in \mathbb{D}.
\end{equation*}

\end{itemize} 
\end{proposition}

\paragraph*{} \textit{Proof of Theorem \ref{GeneFMI}}. Using arguments of \cite[Proposition 2.5]{ArLM}, one can extend (\ref{polbound}), to functions $h$ in $H^\infty_0(\mathbb{D}^d)$. Then we have for any function $g$ in $H^\infty_0(\mathbb{D}^d)$
\begin{equation} \label{HinfiniD}
\left\| g(T_1,...,T_d) \right\| \lesssim \left\| g \right\|_{\infty,\mathbb{D}^d}
\end{equation}
and as the $H^\infty$ functional calculus passes to the adjoint, we also have
\begin{equation} \label{HinfiniD*}
\left\| g(T_1^*,...,T_d^*) \right\| \lesssim \left\| g \right\|_{\infty,\mathbb{D}^d}.
\end{equation}

Let $\gamma_1,...,\gamma_d$ be angles of $(0,\frac{\pi}{2})$ such that every $\gamma_k$ is strictly larger than the $R$-type of $T_k$. Proposition \ref{DecompFMI} yields sequences of functions $(\phi_{i_k,k})_{i_k \geq 1}$, $(\varphi_{i_k,k})_{i_k \geq 1}$, $(\theta_{i_k,k})_{i_k \geq 1}$ and $(\psi_{i_k,k})_{i_k \geq 1}$ of $H^\infty_0(\mathbb{D})$, $k=1,...,d$, such that

\begin{itemize}
\item[(i)] For every $i_k \geq 1$ and $k=1,...,d$, we have 
\begin{equation} \label{fact}
\phi_{i_k,k} = \theta_{i_k,k} \varphi_{i_k,k} \psi_{i_k,k}
\end{equation}
\item[(ii)] There exists a constant $c >0$ such that for $k=1,...,d$ and every $z_k$ in $\mathbb{D}$ 
\begin{equation} \label{est1}
\sum_{i_k=1}^\infty \left| \varphi_{i_k,k}(z_k) \right| \leq c, \qquad \sum_{i_k=1}^\infty \left| \psi_{i_k,k}(z_k) \right| \leq c;
\end{equation}
\item[(iii)] There exists a constant $e >0$ such that for every $i \geq 1$ and $k=1,...,d$
\begin{equation} \label{est2}
\int_{\partial B_{\gamma_k}} \frac{\left| \theta_{i_k,k}(z) \right|}{\left| 1 - z_k \right|} \left| dz_k \right| \leq e;
\end{equation}
\item[(iv)] For $k=1,...,d$ and every $z_k$ in $\mathbb{D}$, the series $\sum_{i_k \geq 1} \phi_{i_k,k}(z_k)$ absolutely converges and there exists a constant $c'$ such that
\begin{equation} \label{est3}
\text{sup} \left\lbrace \sum_{i_k=1}^{\infty} \left| \phi_{i_k,k}(z_k) \right| : z_k \in \mathbb{D} \right\rbrace \leq c'.
\end{equation}
Moreover, we have
\begin{equation} \label{ident}
\sum_{i_k=1}^{\infty} \phi_{i_k,k}(z_k) = 1, \qquad z_k \in \mathbb{D}, \qquad k=1,...,d.
\end{equation}
\end{itemize} 

Let now $b_1,...,b_d$ be angles such that $b_k \in (\gamma_k,\frac{\pi}{2})$, $k=1,...,d$. For the rest of the proof, we fix $h$ a function of $H_0^\infty(B_{b_1} \times \cdots \times B_{b_d})$. By (\ref{ident}), we have the identity 
\begin{equation*} \label{decompseries}
h(z_1,...,z_d) = \sum_{i_1,...,i_d=1}^\infty h(z_1,...,z_d) \phi_{i_1,1}(z_1) \cdots \phi_{i_d,d}(z_d), \qquad (z_1,...,z_d) \in B_{b_1} \times \cdots B_{b_d},
\end{equation*}
where the series of the right hand side is absolutely summable using (\ref{est3}).

Let us now prove that we have 
\begin{equation} \label{decomph(T)}
h(T_1,...,T_d) = \sum_{i_1,...,i_d=1}^\infty h(T_1,...,T_d) \phi_{i_1,1}(T_1) \cdots \phi_{i_d,d}(T_d)
\end{equation}
with absolute convergence of the series. Indeed, write 
\begin{equation*}
h(T_1,...,T_d) = \left( \frac{1}{2\text{i}\pi} \right)^d \int_{\prod_{k=1}^d \partial B_{\gamma_k}} \sum_{i_1,...,i_d=1}^\infty h(z_1,...,z_d) \phi_{i_1,1}(z_1) \cdots \phi_{i_d,d}(z_d) \prod_{k=1}^d R(z_k,T_k) \prod dz_k
\end{equation*}
Hypothesis (\ref{est3}) allows us to write 
\begin{align*}
h(T_1,...,T_d) & =   \sum_{i_1,...,i_d=1}^\infty \left( \frac{1}{2\text{i}\pi} \right)^d \int_{\prod_{k=1}^d \partial B_{\gamma_k}} h(z_1,...,z_d) \phi_{i_1,1}(z_1) \cdots \phi_{i_d,d}(z_d) \prod_{k=1}^d R(z_k,T_k) \prod dz_k \\
 			& = \sum_{i_1,...,i_d=1}^\infty h(T_1,...,T_d) \phi_{i_1,1}(T_1) \cdots \phi_{i_d,d}(T_d)
\end{align*}
and we obtain (\ref{decomph(T)}).

Next we take $x$ in $X$ and $y$ in $X^*$. For convenience, we may write $T=(T_1,...,T_d)$. We compute using (\ref{fact}) and (\ref{decomph(T)})

\begin{align*}
\left\langle h(T)x,y \right\rangle & = \sum_{i_1,...,i_d=1}^\infty \left\langle  h(T) \prod_{k=1}^d \theta_{i_k,k}(T_k) \varphi_{i_k,k}(T_k) \psi_{i_k,k}(T_k) x,y \right\rangle \\
 	  							& = \sum_{i_1,...,i_d=1}^\infty \left\langle  h(T) \prod_{k=1}^d \theta_{i_k,k}(T_k) \varphi_{i_k,k}(T_k)  x, \prod_{k=1}^d \psi_{i_k,k}(T_k)^* y \right\rangle.
\end{align*}

Applying Cauchy-Scwharz inequality, we obtain
\begin{align*}
\left| \left\langle h(T)x,y \right\rangle \right|  \leq & \left\| \sum_{i_1,...,i_d=1}^\infty r_{i_1,...,i_d} \otimes \left[ h(T) \prod_{k=1}^d \theta_{i_k,k}(T_k) \right] \left[ \prod_{k=1}^d \varphi_{i_k,k}(T_k) \right] x\right\|_{\text{Rad}(X)} \\
  & \left\| \sum_{i_1,...,i_d=1}^\infty r_{i_1,...,i_d} \otimes \prod_{k=1}^d \psi_{i_k,k}(T_k)^* y\right\|_{\text{Rad}(X^*)}.
\end{align*}

The sequel is the same as for Theorem \ref{th21} of this paper. Recall that we have
\begin{equation*}
h(T) \prod_{k=1}^d \theta_{i_k,k}(T_k) = \left( \frac{1}{2\text{i}\pi} \right)^d \int_{\prod_{k=1}^d \partial B_{\gamma_k}} h(z_1,...,z_d) \prod_{k=1}^d \theta_{k,i_k}(z_k) R(z_k,T_k) \prod_{k=1}^d dz_k.
\end{equation*}

Using that the set $E$ defined in (\ref{ensRborne}) is $R$-bounded together with assumption (\ref{est2}), we have
\begin{align} \label{hTxy}
\left| \left\langle h(T)x,y \right\rangle \right| \lesssim \left\| h \right\|_{\infty, \prod_{k=1}^d B_{\gamma_k}} & \left\| \sum_{i_1,...,i_d=1}^\infty r_{i_1,...,i_d} \otimes \prod_{k=1}^d \varphi_{i_k,k}(T_k) x\right\|_{\text{Rad}(X)} \times \\
\nonumber                                   &  \left\| \sum_{i_1,...,i_d=1}^\infty r_{i_1,...,i_d} \otimes \prod_{k=1}^d \psi_{i_k,k}(T_k)^* y\right\|_{\text{Rad}(X^*)}.
\end{align}

To conclude, using (\ref{est1}), we claim that for any numbers $s_{i_1,...,i_d}$ in $\left\lbrace -1,1 \right\rbrace$, $(i_1,...,i_d) \in (\mathbb{N}^*)^d$, the sums
\begin{equation*}
S_1(z_1,...,z_d) = \sum_{i_1,...,i_d}^\infty s_{i_1,...,i_d} \prod_{k=1}^d \varphi_{i_k,k}(z_k), \qquad S_2(z_1,...,z_d)=\sum_{i_1,...,i_d}^\infty s_{i_1,...,i_d} \prod_{k=1}^d \psi_{i_k,k}(z_k)
\end{equation*}
are uniformly bounded for $(z_1,...,z_d)$ in $\mathbb{D}^d$. 

Now we can use (\ref{HinfiniD}) with $g=S_1$ and (\ref{HinfiniD*}) with $g=S_2$ to see that taking Rademacher averages we find
\begin{equation*}
\left\| \sum_{i_1,...,i_d=1}^\infty r_{i_1,...,i_d} \otimes \prod_{k=1}^d \varphi_{i_k,k}(T_k) x\right\|_{\text{Rad}(X)} \lesssim \left\| x \right\|
\end{equation*}
and
\begin{equation*}
 \left\| \sum_{i_1,...,i_d=1}^\infty r_{i_1,...,i_d} \otimes \prod_{k=1}^d \psi_{i_k,k}(T_k)^* y \right\|_{\text{Rad}(X^*)} \lesssim \left\| y \right\|.
\end{equation*}

Finally, it follows from all the inequalities above and (\ref{hTxy}) that we have an estimate of type
$$ \left| \left\langle h(T)x,y \right\rangle \right| \lesssim \left\| h \right\|_{\infty, \prod_{k=1}^d B_{\gamma_k}} \left\| x \right\| \left\| y \right\|$$
and $\left\| h(T) \right\| \lesssim \left\| h \right\|_{\infty, \prod_{k=1}^d B_{\gamma_k}} \lesssim \left\| h \right\|_{\infty, \prod_{k=1}^d B_{b_k}}$ is a straightforward consequence. This means that $(T_1,...,T_d)$ has a $H^\infty(B_{b_1} \times \cdots \times B_{b_d})$ joint functional calculus. $\square$

\paragraph*{} \textit{Proof of Theorem \ref{Diljointcalc}}. If (\ref{dilationisom}) is verified, then for any polynomial function $P$ of $d$ variables we have $P(T_1,...,T_d) = Q P(U_1,...,U_d) J$ and then
\begin{equation*}
\left\| P(T_1,...,T_d) \right\| \leq \left\| Q \right\| \left\| J \right\| \left\| P(U_1,...,U_d) \right\|.
\end{equation*}

Now, if $(U_1,...,U_d)$ admits a $C(\mathbb{T}^d)$ bounded calculus, then for any polynomial function $P$ as above we have
\begin{equation*}
\left\| P(U_1,...,U_d) \right\| \lesssim \text{sup} \left\lbrace \left| P(z_1,...,z_d) \right| : (z_1,...,z_d) \in \mathbb{T}^d \right\rbrace.
\end{equation*}

Thus, $(T_1,...,T_d)$ is polynomially bounded and as every $T_k$ is a $R$-Ritt operator, Theorem \ref{GeneFMI} implies that $(T_1,...,T_d)$ admits a $H^\infty(B_{b_1} \times \cdots \times B_{b_d})$ joint functional calculus for some $b_1,...,b_d$ in $(0,\frac{\pi}{2})$. $\square$

\bigskip
\noindent
{\bf Acknowledgements.} 
The author was supported by the French 
``Investissements d'Avenir" program, 
project ISITE-BFC (contract ANR-15-IDEX-03).


\begin{thebibliography}{99}
\bibitem{Al} D. Albrecht, {\it Functional calculi of commuting unbounded operators},
Ph.D. Thesis (Monash University, Melbourne, Australia, 1994).
\bibitem{An} T. Ando, {\it On a pair of commuting contractions}, Acta Sci. Math. 24 (1963), 88-90.
\bibitem{AFLM} C. Arhancet, S. Fackler and C. Le Merdy, 
{\it Isometric dilations and $H^\infty$ calculus for bounded 
analytic semigroups and Ritt operators},  Trans. Amer. Math. Soc. 369 (2017), no. 10, 6899-6933.
\bibitem{AL} C. Arhancet and C. Le Merdy, {\it Dilation of Ritt operators on $L_p$-spaces}, 
Isra\"el J. Math. 201 (2014), no. 1, 373-414.
\bibitem{ArLM} O. Arrigoni and C. Le Merdy, {\it $H^\infty$-functional calculus for commuting families of Ritt operators and sectorial operators},  Oper. Matrices 13 (2019), no. 4, 1055–1090.
\bibitem{ArLM2} O. Arrigoni and C. Le Merdy, {\it New properties of the multivariables $H^\infty$ functional
calculus of sectorial operators}, arXiv 2007.04580v1 [math.FA].
\bibitem{Bu} D. L. Burkholder, {\it Martingales and singular integrals in Banach spaces},
pp. 233-269 in ``Handbook of the geometry of Banach spaces, Vol. I", North-Holland, Amsterdam, 2001.
\bibitem{Fro} A. Fr\"ohlich, L. Weis, {\it $H^\infty$ calculus and dilations}, Bull. Soc. Math. France, \textbf{134} (4), 2006, p 487-508.
\bibitem{FMI} E. Franks and A. McIntosh, {\it Discrete quadratic estimates and 
holomorphic functional calculi in Banach spaces}, Bull. Austral. 
Math. Soc. 58 (1998), 271-290.
\bibitem{GT} A.Gomilko and Y.Tomilov, {\it On discrete subordination of power bounded and Ritt operators},  Indiana Univ. Math. J. 67 (2018), no. 2, 781–829. (Reviewer: Florian Horia Vasilescu) 47A60 (26A48 46L10 47A35 47D03).
\bibitem{Hyt1} T. Hyt\" onen, J. van Neerven, M. Veraar and L. Weis, 
{\it Analysis in Banach spaces I}, Springer, 2016.
\bibitem{Hyt2} T. Hyt\" onen, J. van Neerven, M. Veraar and L. Weis, 
{\it Analysis in Banach spaces II}, Springer, 2016.
\bibitem{KaW} C. Kaiser and L.W. Weis, {\it Wavelet transform for functions with values in UMD space}, Studia Math. 186 (2008), no. 2, 101–126.
\bibitem{K} S. Kwapie{\'n}, {\it On Banach spaces containing $c_0$}, Studia Math. 52, p 187-188, 1974.
\bibitem{KW1} N. J. Kalton and L. Weis,
{\it The $H^\infty$-calculus and sums of closed operators}, Math. Ann. 321 (2001), no. 2, 319-345. 
\bibitem{LLLM} F. Lancien, G. Lancien and C. Le Merdy, {\it A joint functional 
calculus for sectorial operators with commuting resolvents}, Proc. 
London Math. Soc. 77 (1998), no. 3, 387-414.
\bibitem{LM} C. Le Merdy, {\it $H^\infty$ functional calculus and square 
function estimates for Ritt operators}, Rev. Mat. Iberoam. 30 (2014), 1149-1190.
\bibitem{Ly} Y. Lyubich, {\it Spectral localization, power boundedness and invariant 
subspaces under Ritt's type condition}, Studia Math. 134 (1999), no. 2, 153-167. 
\bibitem{Ma} B. Maurey, {\it Type, cotype and $K$-convexity}, Handbook of the geometry of Banach spaces, Vol. 2, 1299–1332, North-Holland, Amsterdam, 2003.
\bibitem{NZ} B. Nagy and J. Zemanek, {\it A resolvent condition implying power boundedness},
Studia Math. 134 (1999), no. 2, 143-151. 
\bibitem{Nev} O. Nevanlinna, {\it Convergence of iterations for linear equations},
Lectures in Mathematics ETH Zürich, Birkhäuser Verlag, Basel, 1993. 
\bibitem{P} G. Pisier, {\it Similarity problems and completely bounded maps}
(Second, expanded edition), Lecture Notes in Mathematics, 1618 Springer-Verlag, 
Berlin, 2001. viii+198 pp. 
\bibitem{P5} G. Pisier, {\it Holomorphic semi-groups and the geometry of Banach spaces}, Annals of Mathematics, no 115 (1982), 375-392.
\end{thebibliography}
\end{document}